\documentclass[12pt,oneside,english,reqno]{amsart}
\usepackage{charter}
\usepackage[T1]{fontenc}
\usepackage[latin9]{inputenc}
\usepackage{geometry}
\usepackage{babel}
\usepackage{prettyref}
\usepackage{mathrsfs}
\usepackage{amsthm}
\usepackage{amstext}
\usepackage{amssymb}
\usepackage{setspace}
\usepackage{esint}
\usepackage[numbers]{natbib}
\usepackage[unicode=true,
 bookmarks=false,
 breaklinks=false,pdfborder={0 0 1},backref=false,colorlinks=false]
 {hyperref}
\hypersetup{pdftitle={Cuntz-Pimsner Algebras for Subproduct Systems},
 pdfauthor={Ami Viselter},
 pdfstartview={XYZ null null 1.25}}

\makeatletter
\theoremstyle{plain}
\newtheorem{thm}{\protect\theoremname}[section]
  \theoremstyle{definition}
  \newtheorem{defn}[thm]{\protect\definitionname}
  \theoremstyle{definition}
  \newtheorem{example}[thm]{\protect\examplename}
  \theoremstyle{plain}
  \newtheorem{prop}[thm]{\protect\propositionname}
  \theoremstyle{plain}
  \newtheorem{lem}[thm]{\protect\lemmaname}
  \theoremstyle{plain}
  \newtheorem{cor}[thm]{\protect\corollaryname}
  \theoremstyle{remark}
  \newtheorem{rem}[thm]{\protect\remarkname}
  \theoremstyle{plain}
  \newtheorem{conjecture}[thm]{\protect\conjecturename}

\usepackage{eucal}
\let\mathcal=\CMcal

\newrefformat{def}{Definition \ref{#1}}
\newrefformat{lem}{Lemma \ref{#1}}
\newrefformat{pro}{Proposition \ref{#1}}
\newrefformat{prop}{Proposition \ref{#1}}
\newrefformat{cor}{Corollary \ref{#1}}
\newrefformat{thm}{Theorem \ref{#1}}
\newrefformat{exa}{Example \ref{#1}}
\newrefformat{rem}{Remark \ref{#1}}
\newrefformat{sec}{\S \ref{#1}}
\newrefformat{conj}{Conjecture \ref{#1}}

\numberwithin{equation}{section}
\DeclareMathOperator{\Img}{Im}
\DeclareMathOperator{\linspan}{span}
\DeclareMathOperator{\clinspan}{\overline{span}}

\newcommand{\br}{{}}

\makeatother

  \providecommand{\conjecturename}{Conjecture}
  \providecommand{\corollaryname}{Corollary}
  \providecommand{\definitionname}{Definition}
  \providecommand{\examplename}{Example}
  \providecommand{\lemmaname}{Lemma}
  \providecommand{\propositionname}{Proposition}
  \providecommand{\remarkname}{Remark}
\providecommand{\theoremname}{Theorem}

\begin{document}
\global\long\def\e{\varepsilon}
\global\long\def\N{\mathbb{N}}
\global\long\def\Z{\mathbb{Z}}
\global\long\def\Q{\mathbb{Q}}
\global\long\def\R{\mathbb{R}}
\global\long\def\C{\mathbb{C}}
\global\long\def\norm#1{\left\Vert #1\right\Vert }

\global\long\def\H{\EuScript H}
\global\long\def\a{\alpha}
\global\long\def\be{\beta}
\global\long\def\l{\lambda}

\global\long\def\tensor{\otimes}

\global\long\def\A{\forall}
\global\long\def\Aa{\mathscr{A}}
\global\long\def\B{\mathscr{B}}

\global\long\def\F{\mathcal{F}}
\global\long\def\L{\mathcal{L}}
\global\long\def\T{\mathcal{T}}
\global\long\def\O{\mathcal{O}}
\global\long\def\K{\mathcal{K}}
\global\long\def\I{\mathcal{I}}
\global\long\def\J{\mathcal{J}}
\global\long\def\P{P}

\global\long\def\moritaeq{\overset{\mathrm{SME}}{\sim}}
\global\long\def\X{\mathsf{X}}
\global\long\def\U{\mathsf{U}}
\global\long\def\M{\mathscr{M}}
\global\long\def\MM{\mathsf{M}}
\global\long\def\leftsub#1#2{\br{}_{#1}#2}

\title{Cuntz-Pimsner Algebras for Subproduct Systems}

\author{Ami Viselter}

\address{Department of Mathematics, Technion -- Israel Institute of Technology,
32000 Haifa, Israel}

\email{viselter@tx.technion.ac.il}
\begin{abstract}
In this paper we generalize the notion of Cuntz-Pimsner algebras of
$C^{*}$-correspondences to the setting of subproduct systems. The
construction is justified in several ways, including the Morita equivalence
of the operator algebras under suitable conditions, and examples are
provided to illustrate its naturality. We also demonstrate why some
features of the Cuntz-Pimsner algebras of $C^{*}$-correspondences
fail to generalize to our setting, and discuss what we have instead.
\end{abstract}
\maketitle
\tableofcontents{}

\section*{Introduction}

The study of the Cuntz-Pimsner algebra of $C^{*}$-correspondences
has its origins in the influential paper of Pimsner \citep{Pimsner}.
Initially defined merely for faithful $C^{*}$-correspondences, the
Cuntz-Pimsner algebra was shown to be a quotient of the Toeplitz algebra
with a special universal property. Katsura \citep{Katsura_2003_construction}
provided a definition for arbitrary $C^{*}$-correspondences, promoting
the Cuntz-Pimsner algebras even further. The construction is flexible
enough to generalize, at the same time, the Cuntz-Krieger algebras,
crossed products by Hilbert $C^{*}$-bimodules \citep{Abadie_Eilers_Exel}
(particularly, crossed products by partial automorphisms) and others.
The Cuntz-Pimsner algebra has been highly popular in research ever
since it was introduced. Many aspects of it have been comprehensively
studied, for instance: $K$-theory \citep{Pimsner,Katsura_2004_On_C_algebras_assoc},
Morita equivalence \citep{Morita_equiv_tensor_alg}, exactness and
nuclearity \citep{Katsura_2004_On_C_algebras_assoc}, ideal structure
\citep{Katsura_2007_ideal_structure_of_C_algebras} and more, and
it has served as a tool in various papers. A crucial step was made
when it was established in \citep{Katsura_2004_On_C_algebras_assoc}
that the Cuntz-Pimsner algebras could be characterized using another
universal property already known for Cuntz-Krieger algebras, namely
the {}``gauge-invariant uniqueness theorem''. This celebrated discovery
(which, particularly, revealed the structure of the isomorphic representations
of the algebra) is so powerful, that it allows an easy proof of many
properties of the Cuntz-Pimsner algebras that have been proven earlier
using more complicated means.

Subproduct systems were introduced in \citep{Subproduct_systems_2009},
where they were studied systematically from several aspects. In our
recent paper \citep{Viselter_cov_rep_subproduct_systems} we continued
these lines, focusing on the tensor and Toeplitz algebras associated
to a subproduct system and their representations. Since these constructions
seem to attract more and more attention recently, it seems that a
generalization of the Cuntz-Pimsner algebra for \emph{subproduct}
systems may have a lot of potential for applications. Nevertheless,
it is not clear \emph{a priori} how the algebra should be defined
so as to have it bear the desirable properties described above.

In this paper we suggest a possible way to define Cuntz-Pimsner algebras
for subproduct systems, which is natural in view of the analysis in,
e.g., \citep[\S 3]{Pimsner}, \citep[\S 4]{Fowler_Muhly_Raeburn_2003}
and \citep[\S 3]{Sims_Yeend_2010}. This has already been done, in
the specific context of a symbolic dynamical system called subshifts,
by Matsumoto (see \citep{Matsumoto_1997} and its follow-ups), and
our definition reduces to his. Throughout the paper we give several
justifications for the {}``correctness'' of our approach. In \prettyref{sec:def_of_algebra}
we explain why the universal characterizations of the Cuntz-Pimsner
algebras for $C^{*}$-correspondences, particularly the gauge-invariant
uniqueness, cannot be employed as-is to the subproduct systems case.
After giving our definition of the algebra, we show that it generalizes
the original construction of Pimsner in two different manners. The
construction is demonstrated by examples in \prettyref{sec:Examples}.
A (partial) characterization of the algebra, in terms of essential
representations of the ambient Toeplitz algebra, is then presented
in \prettyref{sec:tame-subproduct-systems}. In \prettyref{sec:Morita-equivalence}
we define a notion of strong Morita equivalence for subproduct systems,
and show that it implies the Morita equivalence of all associated
operator algebras. This can be seen as another strength of the proposed
definition. The last section is devoted to some open questions.

\section{Preliminaries}

We start with a brief summary of the definitions we need (see \citep{Viselter_cov_rep_subproduct_systems}
and the references therein for a more thorough background), notation
and general assumptions. The reader should be familiar with the basics
of Hilbert $C^{*}$-modules found in \citep[Ch.~1-4]{Lance}. The
notation $\left\langle \cdot,\cdot\right\rangle $ is reserved for
the rigging in Hilbert modules.
\begin{defn}[{\citep[Definition 2.1]{Tensor_algebras}}]
Let $\M,\mathscr{N}$ be $C^{*}$-algebras. A (right) Hilbert $C^{*}$-module
$E$ over $\M$ is called an $\mathscr{N}-\M$ ($C^{*}$-)\emph{ correspondence}
if it is also equipped with a left $\mathscr{N}$-module structure,
implemented by a $*$-homomorphism $\varphi:\mathscr{N}\to\L(E)$;
that is, $a\cdot\zeta:=\varphi(a)\zeta$ for $a\in\mathscr{N}$, $\zeta\in E$.
We say that $E$ is \emph{faithful} if $\varphi$ is faithful and
\emph{essential} if $\varphi(\mathscr{N})E$ is total in $E$. If
$\mathscr{N}=\M$, we say that $E$ is a ($C^{*}$-) correspondence
over $\M$.\end{defn}
\begin{example}
Every Hilbert space is a $C^{*}$-correspondence over $\C$.
\end{example}

\begin{example}
\label{exa:alpha_M}If $\M$ is a $C^{*}$-algebra and $\a$ is an
endomorphism of $\M$, we write $_{\a}\M$ for the $C^{*}$-correspondence
that is equal to $\M$ as sets, with the obvious right $\M$-module
structure, and left $\M$-action given by $\varphi(a)b:=\a(a)b$ for
$a,b\in\M$.
\end{example}

\begin{example}
\label{exa:quiver_corrspndnc}Every quiver (directed graph) possesses
an associated $C^{*}$-correspondence. See \citep[p.~193, (2)]{Pimsner}
or \citep[Example 2.9]{Tensor_algebras} for details.
\end{example}
The full Fock space of a $C^{*}$-correspondence $E$ over $\M$ is
the $C^{*}$-correspondence $\F_{E}:=\bigoplus_{n\in\Z_{+}}E^{\tensor n}$
(also over $\M)$, where $E^{\tensor0}$ equals $\M$ by definition.

The original definition of subproduct systems (\citep[Definitions 1.1, 6.2]{Subproduct_systems_2009})
is for the context of von Neumann algebras. We require its adaptation
to the $C^{*}$-setting of \citep{Viselter_cov_rep_subproduct_systems}.
\begin{defn}
\label{def:subpr_sys}A family $X=\left(X(n)\right)_{n\in\Z_{+}}$
of $C^{*}$-correspondences over a $C^{*}$-algebra $\M$ is called
a (standard) \emph{subproduct system} if $X(0)=\M$ and for all $n,m\in\Z_{+}$,
$X(n+m)$ is an orthogonally-complementable sub-correspondence of
$X(n)\tensor X(m)$. This implies, in particular, that $X(n)$ is
essential for each $n\in\N$.\end{defn}
\begin{example}
If $E$ is an essential $C^{*}$-correspondence over $\M$, the \emph{product}
system $X_{E}$, defined by $X_{E}(n):=E^{\tensor n}$ for each $n\in\Z_{+}$,
is trivially a subproduct system.
\end{example}

\begin{example}[{\citep[Example 1.3]{Subproduct_systems_2009}}]
\label{exa:gen_SSP}Fix a Hilbert space $\H$, and let $X(n):=\H^{\circledS n}$
(the $n$-fold symmetric tensor product of $\H$) for every $n$.
The resulting family $X$ satisfies the requirements of \prettyref{def:subpr_sys}.
It is called the \emph{symmetric subproduct system} over $\H$, and
denoted by $\mathrm{SSP}_{\H}$. Specifically, we put $\mathrm{SSP}_{d}:=\mathrm{SSP}_{\C^{d}}$
for $d\in\N$ and $\mathrm{SSP}_{\infty}:=\mathrm{SSP}_{\ell_{2}(\N)}$.
\end{example}
The reader is urged to consult \citep{Subproduct_systems_2009} for
many other interesting examples of subproduct systems.

Given a subproduct system $X$, we shall use the following notation
throughout the paper. Set $E:=X(1)$. The \emph{$X$-Fock} space is
the sub-correspondence
\[
\F_{X}:=\bigoplus_{n\in\Z_{+}}X(n)
\]
of the full Fock space $\F_{E}$. For all $n\in\Z_{+}$ we have $X(n)\subseteq E^{\tensor n}$.
Let $p_{n}\in\L(E^{\tensor n})$ stand for the (orthogonal) projection
of $E^{\tensor n}$ onto $X(n)$, denote by $Q_{n}\in\L(\F_{X})$
the projection of $\F_{X}$ onto the direct summand $X(n)$, and define
$R_{n}:=Q_{0}+Q_{1}+\ldots+Q_{n}$, $R_{n}':=\bigvee_{k\geq n}Q_{k}$.

The \emph{$X$-shifts} are the operators $S_{n}(\zeta)\in\L(\F_{X})$
($n\in\Z_{+}$, $\zeta\in X(n)$) given by
\[
S_{n}(\zeta)\eta:=p_{n+m}(\zeta\tensor\eta)
\]
for $m\in\Z_{+}$, $\eta\in X(m)$. We write $\varphi_{\infty}(\cdot)$
for $S_{0}(\cdot)$. In case the context is not clear, we will add
the subproduct system letter as a superscript, e.g. $p_{n}^{X}$,
$Q_{n}^{X}$, $S_{n}^{X}$, etc. A direct calculation shows that the
adjoint $S_{n}^{X}(\zeta)^{*}$ is a restriction of the adjoint in
the full Fock space, $S_{n}^{X_{E}}(\zeta)^{*}$, to $\F_{X}$. It
satisfies $S_{n}^{X_{E}}(\zeta)^{*}(\eta_{1}\tensor\eta_{2})=\left\langle \zeta,\eta_{1}\right\rangle \eta_{2}$
for each $\zeta,\eta_{1}\in E^{\tensor n}$ and $\eta_{2}\in E^{\tensor m}$.

The \emph{tensor algebra} $\T_{+}(X)$ is the non-selfadjoint subalgebra
of $\L(\F_{X})$ generated by all $X$-shifts. The \emph{Toeplitz
algebra} $\T(X)$ is the $C^{*}$-subalgebra of $\L(\F_{X})$ generated
by the same operators. It admits a natural action of $\mathbb{T}$,
called the \emph{gauge action}, defined by $\a_{\l}(S_{n}(\zeta)):=\l^{n}S_{n}(\zeta)$
for all $\l\in\mathbb{T}$, $n\in\Z_{+}$ and $\zeta\in X(n)$. It
is useful that the $k$th spectral subset of $\a$ (\citep[Definition 2.1]{Exel_Circle_Actions})
equals the closed linear span of all monomials in $\T(X)$ of degree
$k$, denoted by $\T_{k}(X)$ (see \citep{Viselter_cov_rep_subproduct_systems},
Definition 4.6 and the text surrounding equation (4.13)).
\begin{defn}[{\citep[Definition 2.11]{Tensor_algebras}}]
\label{def:cov_rep_corres}Let $E$ be a $C^{*}$-correspondence
over $\M$ and let $\H$ be a Hilbert space. A pair $(T,\sigma)$
is called a \emph{covariant representation} of $E$ on $\H$ if:
\begin{enumerate}
\item $\sigma$ is a $C^{*}$-representation of $\M$ on $\H$;
\item $T$ is a linear mapping from $E$ to $B(\H)$;
\item and $T$ is a bimodule map with respect to $\sigma$, that is, $T(\zeta a)=T(\zeta)\sigma(a)$
and $T(a\zeta)=\sigma(a)T(\zeta)$ for all $\zeta\in E$ and $a\in\M$.
\end{enumerate}
When this holds, the formula $\widetilde{T}(\zeta\tensor h):=T(\zeta)h$
($\zeta\in E$, $h\in\H)$ defines a contraction $\widetilde{T}:E\tensor_{\sigma}\H\to\H$.
\end{defn}

\begin{defn}[{\citep[Definition 1.5]{Subproduct_systems_2009}, \citep[Definitions 2.6, 2.20]{Viselter_cov_rep_subproduct_systems}}]
Let $X=\left(X(n)\right)_{n\in\Z_{+}}$ be a subproduct system. A
family $T=\left(T_{n}\right)_{n\in\Z_{+}}$ is called a \emph{covariant
representation} of $X$ on $\H$ if:
\begin{enumerate}
\item writing $\sigma:=T_{0}$, the pair $(T_{n},\sigma)$ is a covariant
representation of $X(n)$ on $\H$ (in the sense of \prettyref{def:cov_rep_corres})
for all $n\in\N$;
\item and for every $n,m\in\Z_{+}$, $\zeta\in X(n)$ and $\eta\in X(m)$,
we have
\[
T_{n+m}(p_{n+m}(\zeta\tensor\eta))=T_{n}(\zeta)T_{m}(\eta).
\]

\end{enumerate}
To such a family we associate the operators $\widetilde{T}_{n}:X(n)\tensor_{\sigma}\H\to\H$
as explained above. We say that $T$ is \emph{pure} if the sequence
$\bigl\{\widetilde{T}_{n}\widetilde{T}_{n}^{*}\bigr\}_{n=1}^{\infty}$
converges to zero in the strong operator topology, and \emph{fully
coisometric} if $\widetilde{T}_{1}$ is coisometric.
\end{defn}
If $\pi$ is a representation of $\T_{+}(X)$ (or of $\T(X)$), setting
$T_{n}(\zeta):=\pi(S_{n}(\zeta))$ yields a covariant representation
of $X$. The opposite direction---namely, determining which covariant
representations arise this way---was the main theme of \citep{Viselter_cov_rep_subproduct_systems}.

\textbf{Standing Hypothesis:} throughout this paper, \emph{all subproduct
systems are assumed faithful} ($X$ is faithful if $X(n)$ is faithful
for all $n\in\N$), and \emph{all representations are assumed nondegenerate}.

\section{\label{sec:def_of_algebra}Construction of the algebra}
\begin{defn}
For a (faithful) subproduct system $X$ we denote by $\J$ the ideal
$\varphi^{-1}(\K(E))$ of $\M$.
\end{defn}
When $X$ is a \emph{product} system $X=X_{E}$ (with $E$ faithful),
$\K(\F_{E}\J)$ is an ideal in $\T(E)$ (actually, it is equal to
$\T(E)\cap\K(\F_{E})$), and the Cuntz-Pimsner algebra $\O(E)$ is
defined in \citep{Pimsner} to be the quotient $\T(E)/\K(\F_{E}\J)$.
\begin{prop}
\label{prop:K_F_X_J}Let $X$ be a subproduct system. For all $a\in\J$
we have $\varphi_{\infty}(a)Q_{0}\in\T(X)$. Moreover, the following
subsets of $\L(\F_{X})$ are equal:
\begin{enumerate}
\item $\K(\F_{X}\J)$
\item the ideal of $\T(X)$ generated by $\varphi_{\infty}(\J)Q_{0}$
\item the ideal of $\L(\F_{X})$ generated by $\varphi_{\infty}(\J)Q_{0}$.
\end{enumerate}
As a result, $\K(\F_{X}\J)$ is an ideal of $\T(X)$.
\end{prop}
The proof is almost identical to that of the corresponding assertion
for product systems (compare \citep[Lemma 2.17]{Tensor_algebras}).
Since the proposition is not essential for the rest of the paper,
we omit the details.

Assume that $E$ is a faithful $C^{*}$-correspondence. Katsura's
gauge-invariant uniqueness theorem \citep[Proposition 7.14]{Katsura_2007_ideal_structure_of_C_algebras}%
\footnote{the term {}``gauge-invariant uniqueness theorem'' refers sometimes
to \citep[Theorem 6.4]{Katsura_2004_On_C_algebras_assoc}; but this
is a consequence of \citep[Proposition 7.14]{Katsura_2007_ideal_structure_of_C_algebras}.%
} asserts that the ideal $\K(\F_{E}\J)$ has the following property:
it is the \emph{largest} ideal of $\T(E)$ such that 1) it is gauge
invariant, and 2) its intersection with $\varphi_{\infty}(\M)$ is
$\left\{ 0\right\} $. If $X$ is a subproduct system, it would be
plausible to define $\O(X)$ to be the quotient of $\T(X)$ by such
an ideal. However, the following example demonstrates that such an
ideal fails to exist even in simple cases.
\begin{example}
\label{exa:no_greatest_gauge_invar_ideal}Consider the symmetric subproduct
system $X=\mathrm{SSP}_{2}$ (see \prettyref{exa:gen_SSP}). Suppose
that there exists a largest ideal $\I$, which does not contain the
unit $I\in\T(X)$, and which is gauge invariant. The ideal $\mathbb{K}$
fulfills these two conditions, so we must have $\mathbb{K}\subseteq\I$.
Since $\T(X)/\mathbb{K}$ is canonically isomorphic to $C(\partial B_{2})$
(see \citep[Theorem 5.7]{Arveson_subalgebras_3}), $\I/\mathbb{K}$
can be identified with an ideal of $C(\partial B_{2})$. There thus
exists a nonempty compact set $M$ such that this ideal equals $\mathcal{V}(M)$,
the set of all elements of $C(\partial B_{2})$ vanishing on $M$.
Direct calculation shows that the gauge action $\tilde{\a}$ on $C(\partial B_{2})$,
induced by the gauge action on $\T(X)$, is given by $(\tilde{\a}_{\lambda}(f))(w)=f(\lambda w)$
(for all $\lambda\in\mathbb{T}$, $f\in C(\partial B_{2})$ and $w\in\partial B_{2}$).

Suppose first that $M$ is not contained in a set of the form
\[
M_{a,b}:=\left\{ (z_{1},z_{2})\in\partial B_{2}:\left|z_{1}\right|=a,\left|z_{2}\right|=b\right\}
\]
(for fixed $a,b\geq0$ with $\left|a\right|^{2}+\left|b\right|^{2}=1$).
Choose $w_{1},w_{2}\in M$ with the property that there are scalars
$\alpha,\beta,\gamma\in\C$ such that the polynomial $f(w):=\alpha\left|z_{1}\right|^{2}+\beta\left|z_{2}\right|^{2}+\gamma$
($w=(z_{1},z_{2})$) satisfies $f(w_{1})\neq0$, $f(w_{2})=0$. Write
\[
f(S):=\alpha S_{1}(e_{1})S_{1}(e_{1})^{*}+\beta S_{1}(e_{2})S_{1}(e_{2})^{*}+\gamma I\in\T(X)
\]
and consider the ideal $\I_{1}:=\left\langle \I\cup\left\{ f(S)\right\} \right\rangle $
of $\T(X)$. Since $\I$ and $\left\{ f(S)\right\} $ are gauge invariant,
so is $\I_{1}$. Moreover, $\I_{1}/\mathbb{K}$ is isomorphic to the
ideal $\left\langle \mathcal{V}(M)\cup\left\{ f\right\} \right\rangle $
of $C(\partial B_{2})$. We have $\left\langle \mathcal{V}(M)\cup\left\{ f\right\} \right\rangle \supsetneqq\mathcal{V}(M)$
(thus $\I_{1}\supsetneqq\I$) because $f(w_{1})\neq0$. In addition,
$g(w_{2})=0$ for all $g\in\left\langle \mathcal{V}(M)\cup\left\{ f\right\} \right\rangle $,
hence $1\notin\left\langle \mathcal{V}(M)\cup\left\{ f\right\} \right\rangle $.
In conclusion, we have constructed a gauge-invariant ideal of $\T(X)$,
of which $I$ is not an element, and which strictly contains $\I$.
This is a contradiction.

Assume that $M$ is contained in $M_{a,b}$ for some $a,b$. Since
$M_{a,b}$ fulfills the desirable conditions, we must have $M=M_{a,b}$,
and by symmetry, necessarily $a=b=\frac{1}{\sqrt{2}}$. Now let $f(w):=\left|z_{1}-z_{2}\right|^{2}$,
$w_{1}:=(\frac{1}{\sqrt{2}},-\frac{1}{\sqrt{2}})$ and $w_{2}:=(\frac{1}{\sqrt{2}},\frac{1}{\sqrt{2}})$,
and continue as above to get a contradiction.
\end{example}
The structure of the Toeplitz algebra of a general subproduct system
$X$ is much more complicated than that of a product system, and defining
the Cuntz-Pimsner algebra of $X$ to be $\T(X)/\K(\F_{X}\J)$ is useless,
as could be seen by examples (cf.~\prettyref{rem:CP_not_with_J}
below). We are thus looking for an alternative.
\begin{lem}
\label{lem:two_ideals_gauge_inva}Let $\mathcal{G}_{1}\trianglelefteq\T(X)$
be gauge invariant.
\begin{enumerate}
\item \label{enu:two_ideals_gauge_inva__1}The set $\bigcup_{k\in\Z}\bigl(\mathcal{G}_{1}\cap\T_{k}(X)\bigr)$
is total in $\mathcal{G}_{1}$.
\item \label{enu:two_ideals_gauge_inva__2}If $\mathcal{G}_{2}\trianglelefteq\T(X)$
and $(\mathcal{G}_{1}\cap\T_{0}(X))_{+}\subseteq\mathcal{G}_{2}$,
then $\mathcal{G}_{1}\subseteq\mathcal{G}_{2}$.
\end{enumerate}
\end{lem}
\begin{proof}
We use the routine methods. Let $\left\{ k_{n}\right\} _{n=1}^{\infty}$
denote Fejér's kernel. For $S\in\T(X)$, write $\sigma_{n}(S):=\frac{1}{2\pi}\int_{t=0}^{2\pi}\a_{\l}(S)k_{n}(\l)\,\mathrm{d}t\in\T(X)$
($\l:=e^{it}$). If $S$ is a monomial of degree $m$, then $\a_{\l}(S)=\l^{m}S$,
and so $\sigma_{n}(S)=\left(\frac{1}{2\pi}\int_{t=0}^{2\pi}\l^{m}k_{n}(\l)\,\mathrm{d}t\right)S\xrightarrow[n\to\infty]{}S$.
Hence $\sigma_{n}(S)\to S$ for every (finite) {}``polynomial''
in $\T(X)$. The estimate $\norm{\sigma_{n}(S)}\leq\norm S$ then
yields that actually $\sigma_{n}(S)\to S$ for all $S\in\T(X)$.

Pick $S\in\mathcal{G}_{1}$. For every $k\in\Z$, $\Phi_{k}(S):=\int_{t=0}^{2\pi}\a_{\l}(S)\l^{-k}\,\mathrm{d}t$
belongs to $\mathcal{G}_{1}\cap\T_{k}(X)$. Since $\sigma_{n}(S)$
is a linear combination of $\left\{ \Phi_{k}(S)\right\} _{k=-n}^{n}$,
\prettyref{enu:two_ideals_gauge_inva__1} is established. Moreover,
$\Phi_{k}(S)^{*}\Phi_{k}(S)\in(\mathcal{G}_{1}\cap\T_{0}(X))_{+}$,
and under the assumptions of \prettyref{enu:two_ideals_gauge_inva__2}
this implies that $\Phi_{k}(S)^{*}\Phi_{k}(S)\in\mathcal{G}_{2}$,
and consequently $\Phi_{k}(S)\in\mathcal{G}_{2}$ (by \citep[Theorem I.5.3]{Davidson_C_algebras_by_example}).
Therefore $\sigma_{n}(S)\in\mathcal{G}_{2}$ for all $n$. As a result,
$S\in\mathcal{G}_{2}$.\end{proof}
\begin{thm}
Let $X$ be a subproduct system. Define an ideal $\I'\trianglelefteq\T(X)$
by
\[
\I':=\left\langle S\in\T_{0}(X):\lim_{n\to\infty}\norm{SQ_{n}}=0\right\rangle
\]
and a subset $\I''\subseteq\T(X)$ by
\[
\I'':=\left\{ S\in\T(X):\lim_{n\to\infty}\norm{SQ_{n}}=0\right\} .
\]
Then $\I'=\I''$, and it is gauge invariant. In particular, $\I''\trianglelefteq\T(X)$.
\end{thm}
In comparison with the product system case, $\I'$ and $\I''$ parallel
$\left\langle \varphi_{\infty}(\J)Q_{0}\right\rangle $ and $\T(E)\cap\K(\F_{E})$
(both are equal to $\K(\F_{E}\J)$), respectively.
\begin{proof}
The set $\I''$ is gauge invariant as $\a_{\l}$ is unitarily implemented
by $W_{\l}\in\L(\F_{X})$ given by $\bigoplus_{n\in\Z_{+}}\zeta_{n}\mapsto\bigoplus_{n\in\Z_{+}}\l^{n}\zeta_{n}$,
which commutes with $Q_{n}$ for all $n$.

Let us prove that $\I''$ is an ideal. It is clearly a linear subspace
and a left ideal. It is also norm-closed. Indeed, let $S\in\overline{\I''}$
and $\e>0$ be given, and choose $T\in\I''$ with $\norm{S-T}<\e$.
In particular, $\norm{SQ_{n}-TQ_{n}}<\e$ for all $n$. Since $n_{0}$
can be produced so as to have $\norm{TQ_{n}}<\e$ for all $n\geq n_{0}$,
we have $\norm{SQ_{n}}<2\e$ for all $n\geq n_{0}$, thus $S\in\I''$.
We next verify that $\I''$ is a right ideal. If $S\in\I''$ and $T\in\T(X)$,
we may assume, having proved that $\I''$ is a closed subspace, that
$T$ is a monomial, say of degree $k\in\Z$. Then $T$ maps $X(n)$
to $X(n+k)$ when $n+k\geq0$, and consequently $\norm{STQ_{n}}=\norm{SQ_{n+k}TQ_{n}}\leq\norm{SQ_{n+k}}\cdot\norm T\to0$
as $n\to\infty$.

In conclusion, $\I''\trianglelefteq\T(X)$, and evidently $\I'\subseteq\I''$.
For the converse, we may use \prettyref{lem:two_ideals_gauge_inva}
with $\mathcal{G}_{1}=\I''$ and $\mathcal{G}_{2}=\I'$, as it is
clear that $\I''\cap\T_{0}(X)=\I'$.\end{proof}
\begin{defn}
Let $X$ be a subproduct system. Denote by $\I$ the gauge-invariant
ideal $\I'=\I''$ of $\T(X)$. The \emph{Cuntz-Pimsner algebra} of
$X$ is defined as $\O(X):=\T(X)/\I$.\end{defn}
\begin{cor}
\label{cor:thm_I}For every $S\in\I$, $\norm{SR_{n}'}\xrightarrow[n\to\infty]{}0$.\end{cor}
\begin{proof}
If $S\in\I\cap\T_{k}(X)$ for some $k\in\Z$ then $\norm{SR_{n}'}=\sup_{m\geq n}\norm{SQ_{m}}\xrightarrow[n\to\infty]{}0$.
The proof is complete using \prettyref{enu:two_ideals_gauge_inva__1}
of \prettyref{lem:two_ideals_gauge_inva} and an approximation argument.
\end{proof}
From \prettyref{prop:K_F_X_J} we clearly obtain $\K(\F_{X}\J)\subseteq\T(X)\cap\K(\F_{X})\subseteq\I$.
For product systems the converse also holds, so our definition generalizes
indeed that of Pimsner (\citep{Pimsner}).
\begin{prop}
If $X$ is a (faithful) \emph{product} system $X_{E}$, then $\I=\K(\F_{E}\J)$,
that is, $\O(X_{E})=\O(E)$.\end{prop}
\begin{proof}
On one hand, $\K(\F_{E}\J)\subseteq\I$. On the other hand, $\I$
is gauge invariant, and if $0\neq a\in\M$ then $\norm{\varphi_{\infty}(a)Q_{n}}=\norm a$
for all $n$, and so $\I\cap\varphi_{\infty}(\M)=\left\{ 0\right\} $.
By Katsura's gauge-invariant uniqueness theorem we have $\I\subseteq\K(\F_{E}\J)$.
\end{proof}
Let $\pi$ be a representation of $\T(X)$ on a Hilbert space $\H$.
Since $\I\trianglelefteq\T(X)$, one can decompose $\pi$ as $\pi_{\I}\oplus\pi_{\T(X)/\I}$,
where $\pi_{\I}$ represents $\T(X)$ on the invariant subspace $\clinspan\mbox{\ensuremath{\pi}(\ensuremath{\I})\ensuremath{\H}}$
and $\pi_{\T(X)/\I}$ on its orthogonal complement. Generally, if
$\H'\subseteq\H$ is an invariant subspace for $\pi$, then the subrepresentation
$\pi'$ on $\H'$ has covariant representation $T'$, with $T_{n}':X(n)\tensor\H'\to\H'$
satisfying $\widetilde{T}_{n}'(\zeta\tensor h)=\widetilde{T}_{n}(\zeta\tensor h)$
and $\widetilde{T}_{n}'^{*}h=\widetilde{T}_{n}h$ for all $n\in\Z_{+}$,
$\zeta\in X(n)$ and $h\in\H'$. Consequently $\widetilde{T}_{n}'\widetilde{T}_{n}'^{*}=(\widetilde{T}_{n}\widetilde{T}_{n}^{*})_{|\H'}$.
\begin{prop}
\label{pro:pi_I_is_pure}Let $X$ be a subproduct system whose fibers
are Hilbert spaces (not necessarily finite dimensional) and $\pi$
a $C^{*}$-representation of $\T(X)$ on a Hilbert space $\H$. Then
the representation $\pi_{\I}$ of $\T(X)$ is pure.\end{prop}
\begin{proof}
Denote by $T=\left(T_{n}\right)_{n\in\Z_{+}},C=\left(C_{n}\right)_{n\in\Z_{+}}$
the covariant representations of $\pi,\pi_{\I}$, respectively. To
verify that $\pi_{\I}$ is pure, it is enough to establish that $(\widetilde{C}_{n}\widetilde{C}_{n}^{*}x,x)\to0$
as $n\to\infty$ for each $x=\pi(S)h$ where $S\in\I$ and $0\neq h\in\H$.

Let $x$ be as above. Given $\e>0$, fix $n$ with $\norm{S^{*}R_{n}'}\leq\e/\norm h$
(see \prettyref{cor:thm_I}). Let $\left(e_{\kappa}\right)_{\kappa\in K}$
be an orthonormal base for $X(n)$. Then from \citep[Lemma 3.5]{Viselter_cov_rep_subproduct_systems},
\begin{multline*}
(\widetilde{C}_{n}\widetilde{C}_{n}^{*}x,x)=(\widetilde{T}_{n}\widetilde{T}_{n}^{*}x,x)=\sum_{\kappa}\bigl(\pi(S_{n}(e_{\kappa})S_{n}(e_{\kappa})^{*})x,x\bigr)\\
=\sum_{\kappa}\bigl(\pi(S^{*}S_{n}(e_{\kappa})S_{n}(e_{\kappa})^{*}S)h,h\bigr)=\sum_{\kappa}\bigl(\pi(S^{*}R_{n}'S_{n}(e_{\kappa})S_{n}(e_{\kappa})^{*}R_{n}'S)h,h\bigr).
\end{multline*}
For every finite subset $F\subseteq K$ we have $\sum_{\kappa\in F}S_{n}(e_{\kappa})S_{n}(e_{\kappa})^{*}\leq I$.
Consequently,
\[
\begin{split}\bigl|(\widetilde{C}_{n}\widetilde{C}_{n}^{*}x,x)\bigr| & =\lim_{\substack{F\subseteq K\\
F\text{ is finite}
}
}\bigl|\bigl(\pi(\sum_{\kappa\in F}S^{*}R_{n}'S_{n}(e_{\kappa})S_{n}(e_{\kappa})^{*}R_{n}'S)h,h\bigr)\bigr|\\
 & \leq\lim_{\substack{F\subseteq K\\
F\text{ is finite}
}
}\norm{\sum_{\kappa\in F}S^{*}R_{n}'S_{n}(e_{\kappa})S_{n}(e_{\kappa})^{*}R_{n}'S}\cdot\norm h^{2}\\
 & \leq\norm{S^{*}R_{n}'}^{2}\cdot\norm h^{2}\leq\e^{2}.
\end{split}
\]
This completes the proof.
\end{proof}
Before giving examples, we show in another way why the definition
of $\O(X)$ as $\T(X)/\I$ makes sense. For $n\in\Z_{+}$, consider
$\L(\oplus_{k=0}^{n}X(k))$ as a subspace of $\L(\F_{X})$, and let
$\mathcal{B}:=\bigcup_{n=0}^{\infty}\L(\oplus_{k=0}^{n}X(k))$. Then
$\mathcal{B}$ is a $*$-algebra, whose closure $\overline{\mathcal{B}}$
is a $C^{*}$-subalgebra of $\L(\F_{X})$. Since, additionally, the
inclusion $\overline{\mathcal{B}}\subseteq\L(\F_{X})$ is nondegenerate,
we may consider the multiplier algebra $M(\overline{\mathcal{B}})$
as a $C^{*}$-subalgebra of $\L(\F_{X})$ in the usual manner. It
is straightforward to check that $\T(X)\subseteq M(\overline{\mathcal{B}})$
as the set of monomials is total in $\T(X)$. Denote by $q$ the quotient
map $M(\overline{\mathcal{B}})\to M(\overline{\mathcal{B}})/\overline{\mathcal{B}}$.
Recall that in case $X$ is the product system $X=X_{E}$ (we are
assuming that $E$ is faithful), Pimsner proved in \citep{Pimsner}
that $\O(E)\cong q(\T(E))$. (As a matter of fact, this was the \emph{original}
definition of $\O(E)$).
\begin{prop}
The ideal $\ker q_{|\T(X)}=\overline{\mathcal{B}}\cap\T(X)$ of $\T(X)$
is equal to $\I$. Equivalently, $\O(X)\cong q(\T(X))$.\end{prop}
\begin{proof}
If $S\in\I\cap\T_{0}(X)$ then for every $\e>0$ there exists some
$n_{0}\in\N$ such that, upon defining $T:=SR_{n_{0}}$, we have $T\in\mathcal{B}$
and $\norm{S-T}\leq\e$. Therefore $\I\cap\T_{0}(X)\subseteq\overline{\mathcal{B}}$,
thus $\I\subseteq\overline{\mathcal{B}}\cap\T(X)$ (because $\overline{\mathcal{B}}\cap\T(X)=\ker q_{|\T(X)}\trianglelefteq\T(X)$).
The converse holds similarly: if $S\in\overline{\mathcal{B}}\cap\T(X)$,
then for all $\e>0$ there is an operator $T\in\mathcal{B}$ such
that $\norm{S-T}\leq\e$, and if $n_{0}\in\N$ is such that $T=TR_{n_{0}}$,
then $\norm{S(I-R_{n_{0}})}\leq2\e$, so that $\norm{SQ_{n}}\leq2\e$
for $n>n_{0}$.
\end{proof}

\section{\label{sec:Examples}Examples}

The next theorem demonstrates certain circumstances under which the
ideal $\I$ may be expressed somewhat more explicitly.
\begin{thm}
\label{thm:Q_n_in_T(X)}Let $X$ be a subproduct system.
\begin{enumerate}
\item \label{enu:Q_n_in_T(X)__1}If $Q_{n}\in\T(X)$ for all $n\in\Z_{+}$,
then $\I=\left\langle Q_{n}:n\in\Z_{+}\right\rangle $.
\item \label{enu:Q_n_in_T(X)__2}If, additionally, $I\in\T(X)$ and $\pi$
is a representation of $\T(X)$ whose associated covariant representation
$T$ satisfies $\widetilde{T}_{n}\widetilde{T}_{n}^{*}=\pi(R_{n}')$
for all $n\in\N$, then $\pi$ admits a Wold decomposition---that
is, it is unitarily equivalent to the direct sum of an induced representation
and a fully-coisometric representation.
\end{enumerate}
\end{thm}
\begin{proof}
\prettyref{enu:Q_n_in_T(X)__1} We clearly have $\left\langle Q_{n}:n\in\Z_{+}\right\rangle \subseteq\I$.
Conversely, suppose that $S\in\I$. From \prettyref{cor:thm_I} one
has $\norm{S-SR_{m}}=\norm{SR_{m+1}'}\to0.$ Since $R_{m}\in\left\langle Q_{n}:n\in\Z_{+}\right\rangle $
for all $m$, we get $S\in\left\langle Q_{n}:n\in\Z_{+}\right\rangle $.

\prettyref{enu:Q_n_in_T(X)__2} If $\pi$ is such a representation
of $\T(X)$ on $\H$, consider its decomposition with respect to the
ideal $\I$, $\pi=\pi_{\I}\oplus\pi_{\T(X)/\I}$, as explained above,
and write $\H':=\clinspan\mbox{\ensuremath{\pi}(\ensuremath{\I})\ensuremath{\H}}$.
Since $\pi_{\T(X)/\I}$ factors through $\T(X)/\I$ by construction,
$I-R_{1}'=Q_{0}\in\I$ and $\widetilde{T}_{1}\widetilde{T}_{1}^{*}=\pi(R_{1}')$,
the representation $\pi_{\T(X)/\I}$ is fully coisometric. Denote
by $C=\left(C_{n}\right)_{n\in\Z_{+}}$ the covariant representation
of $\pi_{\I}$. Then $\widetilde{C}_{n}\widetilde{C}_{n}^{*}=\pi(R_{n}')_{|\H'}$,
and in the terminology of \citep[Definition 2.8]{Viselter_cov_rep_subproduct_systems}
we have that $\Delta_{*}(C)=\pi(Q_{0})_{|\H'}$ and
\begin{multline*}
\Delta_{*}(C)C_{n}(\zeta)^{*}C_{n}(\zeta)\Delta_{*}(C)=\pi(Q_{0}S_{n}(\zeta)^{*}S_{n}(\zeta)Q_{0})_{|\H'}\\
=\pi(\varphi_{\infty}(\left\langle \zeta,\zeta\right\rangle )Q_{0})_{|\H'}=C_{0}(\left\langle \zeta,\zeta\right\rangle )\Delta_{*}(C).
\end{multline*}
In other words, $C$ is relatively isometric (\citep[Definition 3.3]{Viselter_cov_rep_subproduct_systems}).
To verify that it is pure, it is enough to establish that $\widetilde{C}_{n}\widetilde{C}_{n}^{*}x\to0$
as $n\to\infty$ for vectors $x$ of the form $\pi(S)h$, $S\in\I$
and $h\in\H$. Indeed, $\widetilde{C}_{n}\widetilde{C}_{n}^{*}x=\pi(R_{n}'S)h\to0$
from \prettyref{cor:thm_I}.

In conclusion, all conditions of \citep[Theorem 3.8]{Viselter_cov_rep_subproduct_systems}
are satisfied, so that $C$ (equivalently, $\pi_{\I}$) is induced.\end{proof}
\begin{cor}
\label{cor:subpr_sys_finite_dim_Hilbert_sp}If $X$ is a subproduct
system of finite dimensional Hilbert spaces, then $\I=\mathbb{K}$
(the compacts over the separable Hilbert space $\F_{X}$), and $X$
fulfills the requirements of \prettyref{thm:Q_n_in_T(X)} for every
representation.\end{cor}
\begin{proof}
From \citep[Proposition 8.1]{Subproduct_systems_2009} it follows
that $\mathbb{K}\subseteq\T(X)$, and it is easily seen that $\left\langle Q_{n}:n\in\Z_{+}\right\rangle =\mathbb{K}$.
The second assertion is a consequence of \citep[Lemma 3.5]{Viselter_cov_rep_subproduct_systems}.\end{proof}
\begin{example}
Take $d\in\N$. By \citep[Theorem 5.7]{Arveson_subalgebras_3} we
get the expected result $\O(\mathrm{SSP}_{d})\cong C(\partial B_{d})$
(see \prettyref{exa:gen_SSP}).
\end{example}

\begin{example}
Let $\Lambda$ be a \emph{subshift} in the sense of \citep{Matsumoto_1997}.
Then the $C^{*}$-algebra $\O_{\Lambda}$ associated with $\Lambda$
is equal to $\O(X_{\Lambda})$, where $X_{\Lambda}$ is the subproduct
system associated with $\Lambda$ as in \citep[\S 12]{Subproduct_systems_2009}
(see Definition 12.1 and Remark 12.2 there).
\end{example}
The conditions of \prettyref{thm:Q_n_in_T(X)} also hold for subproduct
systems whose fibers are \emph{not} Hilbert spaces. \prettyref{exa:positive_P}
below is an illustration of this.

Subproduct systems whose fibers are \emph{infinite} dimensional Hilbert
spaces, which do not satisfy the conditions of \prettyref{thm:Q_n_in_T(X)},
are also of interest. We next consider the subproduct system $\mathrm{SSP}_{\infty}$.
The following lemma is required to express $\I$ concretely.
\begin{lem}
\label{lem:Abelian_C_algebra_gen_by_shifts}Let $\Aa$ be a unital
Abelian $C^{*}$-algebra. Suppose that there exists a bounded linear
mapping $A:\ell_{2}(\N)\to\Aa$ such that:
\begin{enumerate}
\item \label{enu:lem:Abelian_C_algebra_gen_by_shifts_1}$\norm{A(e_{n})}=1$
for each $n\in\N$
\item \label{enu:lem:Abelian_C_algebra_gen_by_shifts_2}$\Aa$ is generated
by $\left\{ I,A(e_{1}),A(e_{2}),\ldots\right\} $
\item \label{enu:lem:Abelian_C_algebra_gen_by_shifts_3}the inequality $A(e_{1})^{*}A(e_{1})+\ldots+A(e_{n})^{*}A(e_{n})\leq I$
holds for all $n\in\N$
\item \label{enu:lem:Abelian_C_algebra_gen_by_shifts_4}for every unitary
$U\in B(\ell_{2}(\N))$, the mapping $A(x)\mapsto A(Ux)$ extends
to an automorphism $\a_{U}$ of $\Aa$.
\end{enumerate}
Then the structure space of $\Aa$ can be naturally identified with
the unit ball
\[
B:=\Bigl\{\left(z_{n}\right)_{n\in\N}\in\overline{\mathbb{D}}^{\N}:\sum_{n=1}^{\infty}\left|z_{n}\right|^{2}\leq1\Bigr\}
\]
of $\ell_{2}$ endowed with the Tychonoff topology.\end{lem}
\begin{proof}
Denote the structure space of $\Aa$ by $M$. From assumptions \ref{enu:lem:Abelian_C_algebra_gen_by_shifts_2}
and \ref{enu:lem:Abelian_C_algebra_gen_by_shifts_3} it follows (see
\citep[Theorem IX.2.11]{DS2}, for example) that the map $\rho\mapsto\left(\rho(A(e_{n}))\right)_{n\in\N}$
is a (topological) embedding of $M$ into $B$. We should prove that
it is surjective.

By \ref{enu:lem:Abelian_C_algebra_gen_by_shifts_1} there is a pure
state (one-dimensional representation) $\rho_{1}$ of $\Aa$ with
$\left|\l\right|=1$ where $\l:=\rho_{1}(A(e_{1}))$. We must therefore
have $\rho_{1}(A(e_{k}))=0$ for all $k\geq2$. Write $B':=\left\{ z\in B:\norm z_{2}=1\right\} $.
Given $z=\left(z_{n}\right)_{n\in\N}\in B'$, let $U\in B(\ell_{2}(\N))$
be a unitary with $(Ue_{n},e_{1})=\overline{\l}z_{n}$ for all $n\in\N$.
The pure state $\rho_{1}\circ\a_{U}$ (see \ref{enu:lem:Abelian_C_algebra_gen_by_shifts_4})
satisfies
\[
(\rho_{1}\circ\a_{U})(A(e_{n}))=\rho_{1}(A(Ue_{n}))=\sum_{k=1}^{\infty}(Ue_{n},e_{k})\rho_{1}(A(e_{k}))=z_{n}\qquad(\A n\in\N),
\]
and consequently $z$ belongs to $M$. Since $B'$ is dense in $B$
and $M$ is closed, we have $B=M$, as desired.\end{proof}
\begin{example}
\label{exa:SSP_infty}Consider the subproduct system $X:=\mathrm{SSP}_{\infty}$.
Its Cuntz-Pimsner algebra is the commutative counterpart of $\O_{\infty}$.
Let $\I_{1}$ denote the ideal in $\T(X)$ generated by the commutators
$\left[S_{1}(e_{n}),S_{1}(e_{m})^{*}\right]$, $n,m\in\N$. As in
\citep[Proposition 5.3]{Arveson_subalgebras_3}, one sees that every
such commutator is in $\I$, so that $\I_{1}\subseteq\I$. The quotient
$\T(X)/\I_{1}$ is a unital Abelian $C^{*}$-algebra.

We would like to apply \prettyref{lem:Abelian_C_algebra_gen_by_shifts}
to $\T(X)/\I_{1}$ and $\T(X)/\I$ with $A(x)$ being defined as $S_{1}(x)+\I_{1}$
and $S_{1}(x)+\I$, respectively. It follows from the definition of
$\I$ that for all $n\in\N$ and $T\in\I$,
\[
\norm{S_{1}(e_{n})+T}\ge\lim_{m\to\infty}\norm{(S_{1}(e_{n})+T)Q_{m}}=\lim_{m\to\infty}\norm{S_{1}(e_{n})Q_{m}}=1
\]
(because $S_{1}(e_{n})(e_{n}^{\tensor m})=e_{n}^{\tensor(m+1)}$ for
all $m$). Therefore $\norm{S_{1}(e_{n})+\I}=1$, and thus $\norm{S_{1}(e_{n})+\I_{1}}=1$,
proving \ref{enu:lem:Abelian_C_algebra_gen_by_shifts_1}. Assumptions
\ref{enu:lem:Abelian_C_algebra_gen_by_shifts_2} and \ref{enu:lem:Abelian_C_algebra_gen_by_shifts_3}
clearly hold in both cases.

To establish \ref{enu:lem:Abelian_C_algebra_gen_by_shifts_4}, let
$U\in B(\ell_{2}(\N))$ be unitary. Define a unitary $W\in B(\F_{X})$
to be the restriction to $\F_{X}$ of the unitary $\bigoplus_{n\in\Z_{+}}U^{\tensor n}$
over the full Fock space $\F_{E}$. The automorphism $\a_{U}$ of
$\T(X)$ mapping $S_{1}(x)$ to $S_{1}(Ux)$ ($x\in E$) is implemented
by $W$. Direct calculation shows that $\a_{U}(\left[S_{1}(e_{n}),S_{1}(e_{m})^{*}\right])\in\I_{1}$
for all $n,m\in\N$, thus $\a_{U}(\I_{1})=\I_{1}$. Furthermore, $\a_{U}(\I)=\I$
as $W$ commutes with $Q_{n}$ for all $n$.

In conclusion, it follows from \prettyref{lem:Abelian_C_algebra_gen_by_shifts}
that $\I_{1}=\I$, and that we have the exact sequence
\[
0\to\I\to\T(X)\to C(B)\to0
\]
(compare \citep[Theorem 5.7]{Arveson_subalgebras_3}). Consequently,
$\O(X)\cong C(B)$.

For a given point $\left(z_{n}\right)_{n\in\N}=z\in B$, the corresponding
representation $\rho_{z}:\T(X)\to\C$ (pulled back from $\T(X)/\I$)
satisfies $\rho_{z}(S_{1}(e_{n}))=z_{n}$ for all $n\in\N$. Therefore,
denoting by $T$ the suitable covariant representation of $\T(X)$,
we obtain $\widetilde{T}_{1}\widetilde{T}_{1}^{*}=\sum_{n=1}^{\infty}\left|z_{n}\right|^{2}=\norm z_{2}^{2}$
(for general $m\in\N$, $\widetilde{T}_{m}\widetilde{T}_{m}^{*}=\sum_{\a\in\N^{m}}|z_{a_{1}}|^{2}\cdots|z_{\a_{n}}|^{2}=\norm z_{2}^{2m}$).
Hence, if $z\neq0$, then $\widetilde{T}_{1}:E\to\C$ is a partial
isometry if and only if $z\in B'$ (see the notation of the proof
of \prettyref{lem:Abelian_C_algebra_gen_by_shifts}), if and only
if $T$ is fully coisometric. In particular, there is an abundance
of representations of $\T(X)$ whose associated covariant representations
are not a partial isometry. In the pathological case $z=0$ we have
$\rho_{z}(S_{1}(e_{n}))=0$ for all $n$ and $\widetilde{T}_{1}=0$.
The covariant representation $T$ is trivially pure, but it is by
no means relatively isometric (see \citep[\S 3]{Viselter_cov_rep_subproduct_systems}).
Particularly, $T$ extends to a $C^{*}$-representation although the
conditions of \citep[Theorem 3.8]{Viselter_cov_rep_subproduct_systems}
are not satisfied. \end{example}
\begin{rem}
\label{rem:CP_not_with_J}The last example shows very clearly that
for subproduct systems, defining the Cuntz-Pimsner algebra as $\T(X)/\K(\F_{X}\J)$
is counter-intuitive, since $\J$ of $\mathrm{SSP}_{\infty}$ is $\left\{ 0\right\} $.
This stands in stark contrast to the Cuntz algebra $\O_{\infty}$,
which equals its corresponding Toeplitz algebra $\T\O_{\infty}$.\end{rem}
\begin{example}[The subproduct system of a {}``positive'' matrix]
\label{exa:positive_P}For a unital $C^{*}$-algebra $\M$ and a
completely positive map $P$ over $\M$, the $C^{*}$-correspondence
$\M\tensor_{P}\M$ over $\M$ (see \citep[\S 5]{Paschke_1973}) is
constructed from the algebraic $\C$-balanced tensor product $\M\tensor_{\textrm{alg}}\M$
by giving it the standard left and right actions and the rigging
\[
\left\langle a\tensor_{P}b,c\tensor_{P}d\right\rangle =b^{*}P(a^{*}c)d.
\]
A sufficient condition for $\M\tensor_{P}\M$ to be faithful is that
$P$ be faithful.

Let $P_{1},P_{2}$ be two completely positive maps over $\M$. Then
$P_{2}P_{1}$ is also a completely positive map over $\M$, and there
is a correspondence isometry
\[
V_{P_{1},P_{2}}:\M\tensor_{P_{2}P_{1}}\M\to\left(\M\tensor_{P_{1}}\M\right)\tensor\left(\M\tensor_{P_{2}}\M\right)
\]
 defined by $a\tensor_{P_{2}P_{1}}b\mapsto\left(a\tensor_{P_{1}}I_{\M}\right)\tensor\left(I_{\M}\tensor_{P_{2}}b\right)$,
$a,b\in\M$.

Henceforth we take $\M:=\C^{d}$, $d\in\N$. In this case, a linear
map $P:\M\to\M$ can be identified with a matrix $P=\left(P_{ij}\right)\in M_{d}(\C)$.
The map $P$ is completely positive if and only if it is positive
(as $\M$ is commutative), which is equivalent to that $P_{ij}\geq0$
for all $i,j$. We also assume that $P$ is faithful, equivalently:
every column of $P$ has at least one entry with value strictly greater
than zero.

Let $e_{1},\ldots,e_{d}$ be the standard basis of $\C^{d}$. Write
$e_{ij}$ for (the equivalence class of) the element $e_{i}\tensor e_{j}$
of $\M\tensor_{P}\M$. Notice that
\begin{equation}
\left\langle e_{ij},e_{kl}\right\rangle _{\M\tensor_{P}\M}=e_{j}^{*}P(e_{i}^{*}e_{k})e_{l}=\begin{cases}
P_{ji}e_{j} & \text{if }(i,j)=(k,l)\\
0 & \text{else.}
\end{cases}\label{eq:M_P_M_rigging}
\end{equation}
In particular, $e_{ij}\neq0$ in $\M\tensor_{P}\M$ if and only if
$P_{ji}>0$.
\end{example}
Let now $G_{\P}$ stand for the quiver with vertices $1,2,\ldots,d$,
and with an edge going from $j$ to $i$ (denoted by $g_{ij})$ if
and only if $P_{ji}>0$. This quiver is the \emph{support} of $\P$.
Write $f_{i}$ for the element of $C(G_{P}^{(0)})\cong\M$ mapping
$i$ to $1$ and all other vertices to $0$, and $f_{ij}$ for the
element of $C(G_{\P}^{(1)})$ mapping $g_{ij}$ to $1$ and all other
edges to $0$ (if $P_{ji}>0$; otherwise, set $f_{ij}:=0$). Then
the $C^{*}$-correspondence $C(G_{\P}^{(1)})$ of $G_{\P}$ (see \prettyref{exa:quiver_corrspndnc})
is naturally isomorphic to $\M\tensor_{\P}\M$ via $\Psi_{\P}:\M\tensor_{\P}\M\to C(G_{\P}^{(1)})$
defined by $e_{ij}\mapsto\sqrt{\P_{ji}}f_{ij}$ .

For all $n\in\N$, the map $\P^{n}$ is (completely) positive over
$\M$ and faithful. Let $\left(P_{ij}^{n}\right)\in M_{d}(\C)$ be
its representing matrix, and denote by $X(n)$ the $C^{*}$-correspondence
$\M\tensor_{\P^{n}}\M$. Write also $X(0):=\M$. Fix $n,m\in\N$.
Regarding $e_{ij}$, $e_{kl}$ as elements of $X(n)$, $X(m)$, respectively,
one sees from \prettyref{eq:M_P_M_rigging} that
\begin{multline*}
\left\langle e_{ij}\tensor e_{kl},e_{ij}\tensor e_{kl}\right\rangle _{X(n)\tensor X(m)}=\left\langle e_{kl},\left\langle e_{ij},e_{ij}\right\rangle _{X(n)}\cdot e_{kl}\right\rangle _{X(m)}\\
=P_{ji}^{n}\left\langle e_{kl},e_{j}\cdot e_{kl}\right\rangle _{X(m)}=P_{lk}^{m}P_{ji}^{n}e_{l}\delta_{j,k}.
\end{multline*}
In particular, $e_{ij}\tensor e_{kl}\neq0$ in $X(n)\tensor X(m)$
if and only if $j=k$ and $P_{ji}^{n},P_{lk}^{m}>0$.

As seen above, we may regard $X(n+m)$ as a sub-correspondence of
$X(n)\tensor X(m)$ via the embedding $V_{n,m}:=V_{\P^{n},\P^{m}}$.
Now $V_{n,m}$ is adjointable, and its adjoint $V_{n,m}^{*}:X(n)\tensor X(m)\to X(n+m)$
is given by
\[
V_{n,m}^{*}(e_{ij}\tensor e_{kl})=\begin{cases}
(P_{li}^{n+m})^{-1}\P_{lj}^{m}\P_{ji}^{n}e_{il} & \text{if }j=k\text{ and }P_{li}^{n+m}>0\\
0 & \text{else.}
\end{cases}
\]
Indeed, for elements of the form $e_{ij}\in X(n)$, $e_{kl}\in X(m)$
and $e_{pq}\in X(n+m)$ we have by \prettyref{eq:M_P_M_rigging}
\begin{multline*}
\left\langle V_{n,m}e_{pq},e_{ij}\tensor e_{kl}\right\rangle _{X(n)\tensor X(m)}=\left\langle \left(e_{p}\tensor_{\P^{n}}I_{\M}\right)\tensor\left(I_{\M}\tensor_{\P^{m}}e_{q}\right),e_{ij}\tensor e_{kl}\right\rangle _{X(n)\tensor X(m)}\\
\begin{split} & =\sum_{t=1}^{d}\left\langle e_{pt}\tensor e_{tq},e_{ij}\tensor e_{kl}\right\rangle _{X(n)\tensor X(m)}=\sum_{t=1}^{d}\left\langle e_{tq},\left\langle e_{pt},e_{ij}\right\rangle _{X(n)}\cdot e_{kl}\right\rangle _{X(m)}\\
 & =P_{ji}^{n}\delta_{p,i}\left\langle e_{jq},e_{j}\cdot e_{kl}\right\rangle _{X(m)}=\delta_{p,i}\delta_{q,l}\delta_{j,k}P_{lj}^{m}P_{ji}^{n}e_{l}.
\end{split}
\end{multline*}
It is easy to check that $(V_{n,m}\tensor I_{X(k)})V_{n+m,k}=(I_{X(n)}\tensor V_{m,k})V_{n,m+k}$
for all $n,m,k$, making $X=(X(n))_{n\in\Z_{+}}$ a subproduct system.

For $n\in\N$, let $Y(n)$ denote the $C^{*}$-correspondence of $G_{P^{n}}$.
Then $Y=\left(Y(n)\right)_{n\in\Z_{+}}$ is a subproduct system with
respect to the embeddings
\[
W_{n,m}:=\left(\Psi_{\P^{n}}\tensor\Psi_{\P^{m}}\right)V_{n,m}\Psi_{\P^{n+m}}^{-1}:Y(n+m)\to Y(n)\tensor Y(m),
\]
which satisfy
\[
W_{n,m}f_{ij}=\frac{1}{\sqrt{\P_{ji}^{n+m}}}\sum_{t=1}^{d}\left(\Psi_{\P^{n}}\tensor\Psi_{\P^{m}}\right)e_{it}\tensor e_{tj}=\sum_{t=1}^{d}\sqrt{\frac{P_{jt}^{m}P_{ti}^{n}}{\P_{ji}^{n+m}}}f_{it}\tensor f_{tj}
\]
and
\[
W_{n,m}^{*}\left(f_{ik}\tensor f_{kl}\right)=\begin{cases}
\sqrt{\frac{P_{lk}^{m}P_{ki}^{n}}{P_{li}^{n+m}}}f_{il} & \text{if }P_{li}^{n+m}>0\\
0 & \text{else.}
\end{cases}
\]

Abbreviate $S_{n}^{Y}(\zeta)$ by $S_{n}(\zeta)$. Let $n,m\in\N$.
Regard $f_{ij}$, $f_{kl}$ as elements of $Y(n)$, $Y(m)$, respectively.
Then
\[
S_{n}(f_{ij})f_{kl}=W_{n,m}^{*}(f_{ij}\tensor f_{kl})=\begin{cases}
\sqrt{\frac{P_{lk}^{m}P_{ki}^{n}}{P_{li}^{n+m}}}\delta_{j,k}f_{il} & P_{li}^{n+m}>0\\
0 & \text{else.}
\end{cases}\in Y(n+m)
\]
Regarding $f_{ij}$, $f_{kl}$ as elements of $Y(n)$, $Y(n+m)$,
respectively, we obtain
\[
S_{n}^{*}(f_{ij})f_{kl}=\sum_{t=1}^{d}\sqrt{\frac{P_{lt}^{m}P_{tk}^{n}}{\P_{lk}^{n+m}}}\left\langle f_{ij},f_{kt}\right\rangle _{Y(n)}\cdot f_{tl}=\sqrt{\frac{P_{lj}^{m}P_{jk}^{n}}{\P_{lk}^{n+m}}}\delta_{i,k}f_{jl}\in Y(m)
\]
while if $f_{ij},f_{kl}\in Y(n)$ then
\[
S_{n}^{*}(f_{ij})f_{kl}=\left\langle f_{ij},f_{kl}\right\rangle _{Y(n)}=\delta_{(i,j),(k,l)}f_{l}.
\]

Given $n,m\in\N$, consider $f_{kl}$ as an element of $Y(n+m)$ (assuming
$P_{lk}^{n+m}>0$). Then
\begin{equation}
\sum_{s,t=1}^{d}S_{n}(f_{ts})S_{n}(f_{ts})^{*}f_{kl}=\sum_{t=1}^{d}\delta_{t,k}\sum_{s=1}^{d}\sqrt{\frac{P_{ls}^{m}P_{sk}^{n}}{\P_{lk}^{n+m}}}S_{n}(f_{ks})f_{sl}=\sum_{s=1}^{d}\frac{P_{ls}^{m}P_{sk}^{n}}{P_{lk}^{n+m}}f_{kl}=f_{kl}.\label{eq:S_n_S_n_star}
\end{equation}
Similarly, it is interesting to note that for all $n\in\N$ the left
multiplication in $Y(n)$ is implemented by compacts: $\varphi(f_{t})=\sum_{s=1}^{d}f_{ts}\tensor f_{ts}^{*}$.
\begin{prop}
\label{prop:positive_matrix_subpr_sys}Let $P\in M_{d}(\C)$ be as
in the last example, and $Y$ be the associated subproduct system.
Let $\pi$ be a representation of $Y$ on a Hilbert space $\H$, with
$T$ the associated covariant representation. Then $\left\{ Q_{n}:n\in\Z_{+}\right\} \subseteq\T(Y)$
and for all $n$,
\begin{equation}
\widetilde{T}_{n}^{*}h=\sum_{i,j=1}^{d}f_{ij}\tensor T_{n}(f_{ij})^{*}h\qquad(\A h\in\H)\label{eq:positive_matrix_subpr_sys_1}
\end{equation}
and
\[
\widetilde{T}_{n}\widetilde{T}_{n}^{*}=\pi(R_{n}').
\]
Hence, $Y$ fulfills the requirements of \prettyref{thm:Q_n_in_T(X)}
for every representation.\end{prop}
\begin{proof}
Equation \prettyref{eq:positive_matrix_subpr_sys_1} is checked by
a simple calculation, because
\[
\left\langle f_{ij},f_{kl}\right\rangle _{Y(n)}=\begin{cases}
\delta_{(i,j),(k,l)}f_{j} & P_{ji}^{n}>0\text{, equivalently: }f_{ij}\neq0\\
0 & \text{else.}
\end{cases}
\]
The other assertions follow from \prettyref{eq:S_n_S_n_star}. We
omit the details.
\end{proof}

\section{\label{sec:tame-subproduct-systems}Essential and fully-coisometric
representations}

As further justification for the definition of $\O(X)$ as $\T(X)/\I$,
we sought  a suitable {}``universality'' property of $\I$, which
could replace the gauge-invariant uniqueness theorem (cf.~\prettyref{exa:no_greatest_gauge_invar_ideal}).
More specifically, our goal was to express $\I$ as the intersection
of a certain set of ideals, as in the next proposition. Unfortunately,
we could generally establish only half of this characterization in
\prettyref{thm:I_subset_cap}. Nevertheless, we exemplify many subproduct
systems for which $\I$ has this property, the most non-standard of
which is the infinite-dimensional symmetric subproduct system $\mathrm{SSP}_{\infty}$
(\prettyref{exa:SSP_infty}).
\begin{prop}
Let $E$ be a faithful and essential $C^{*}$-correspondence. Then
the intersection of the kernels of all fully-coisometric $C^{*}$-representations
of $\T(E)$ is $\K(\F_{E}\J)$.\end{prop}
\begin{proof}
In case $E$ is \emph{full}, this result is a reformulation of \citep[Theorem 1.2]{Hirshberg_2005_Ess_Rep}.

For the general case, denote the above-mentioned intersection by $\mathcal{P}$.
Let $\pi$ be a fully-coisometric $C^{*}$-representation of $\T(E)$.
Then $\K(\F_{E}\J)\subseteq\ker\pi$ by \citep[Lemma 5.5]{Tensor_algebras}.
Moreover, if $\l\in\mathbb{T}$, then the $C^{*}$-representation
$\pi\circ\a_{\l}$ of $\T(E)$ is also fully coisometric. As a result,
$\mathcal{P}$ is gauge invariant. By \citep[Theorem 8.3]{Skeide_2009_Unit_vec_Morita_equi_end},
there exists a fully-coisometric $C^{*}$-representation $\pi$ of
$\T(E)$ such that $\pi\circ\varphi_{\infty}$ is faithful. Consequently,
$\mathcal{P}\cap\varphi_{\infty}(\M)=\left\{ 0\right\} $. Katsura's
gauge-invariant uniqueness theorem therefore implies that $\mathcal{P}=\K(\F_{E}\J)$,
as desired.\end{proof}
\begin{defn}
Let $X$ be a subproduct system. A $C^{*}$-representation $\pi$
of $\T(X)$ on $\H$ is said to be \emph{essential} if the associated
covariant representation $T$ satisfies that $\Img\widetilde{T}_{n}$
is dense in $\H$ (equivalently: $\bigcup_{\zeta\in X(n)}\Img T_{n}(\zeta)$
is total in $\H$) for all $n$.
\end{defn}
This requirement is weaker than $\pi$ being fully coisometric, and
it is often strictly weaker; see \prettyref{exa:SSP_infty} (also
compare \citep[Remark 4.2]{Viselter_cov_rep_subproduct_systems}).
Nevertheless, in some special cases, $\pi$ is essential if and only
if it is fully coisometric. This happens, in particular, when the
operators $\widetilde{T}_{n}$ are automatically partial isometries.
For instance:
\begin{enumerate}
\item if $X$ is a \emph{product} system, because then the operators $\widetilde{T}_{n}$
are isometries;
\item if $\pi$ satisfies the conditions of \prettyref{thm:Q_n_in_T(X)},
\prettyref{enu:Q_n_in_T(X)__2}; for example, if the fibers of $X$
are finite-dimensional Hilbert spaces, or if $X$ is the subproduct
system over $\C^{d}$ constructed in \prettyref{exa:positive_P} (by
\prettyref{prop:positive_matrix_subpr_sys}).\end{enumerate}
\begin{thm}
\label{thm:I_subset_cap}If $X$ is a subproduct system and $\pi$
is an essential $C^{*}$-representation of $\T(X)$, then $\I\subseteq\ker\pi$.\end{thm}
\begin{proof}
Suppose that $\pi$ represents $\T(X)$ on the Hilbert space $\H$.
Let $S\in\I$. We have to show that $S\in\ker\pi$. Fix $x\in\H$
and $\e>0$, and choose $n$ such that $\norm{SR_{n}'}\leq\e$ (see
\prettyref{cor:thm_I}). By assumption, the set $\linspan\left\{ \pi(S_{n}(\zeta))y:\zeta\in X(n),y\in\H\right\} $
is dense in $\H$, so there exist $t\in\N$, $\zeta_{1},\ldots,\zeta_{t}\in X(n)$
and $y_{1},\ldots y_{t}\in\H$ so that
\[
\norm{x-z}_{\H}\leq\e\text{ for }z:=\pi(S_{n}(\zeta_{1}))y_{1}+\ldots+\pi(S_{n}(\zeta_{t}))y_{t}.
\]
Then
\[
\begin{split}\norm{\pi(S)z}_{\H}^{2} & =\sum_{i,j=1}^{t}\bigl(\pi(S_{n}(\zeta_{j})^{*}S^{*}SS_{n}(\zeta_{i}))y_{i},y_{j}\bigr)_{\H}\\
 & =\left(\pi^{(t)}\left(\bigl(S_{n}(\zeta_{i})^{*}S^{*}SS_{n}(\zeta_{j})\bigr)_{i,j=1}^{t}\right)\left(y_{k}\right)_{k=1}^{t},\left(y_{\ell}\right)_{\ell=1}^{t}\right)_{\H\tensor\C^{t}}.
\end{split}
\]
For all $\zeta\in X(n)$ we have $SS_{n}(\zeta)=SR_{n}'S_{n}(\zeta)$
and $0\leq R_{n}'S^{*}SR_{n}'\leq\e^{2}I_{\F_{X}}$. Hence, using
the (positive) matrix inequality
\[
\bigl(S_{n}(\zeta_{i})^{*}R_{n}'S^{*}SR_{n}'S_{n}(\zeta_{j})\bigr)_{i,j=1}^{t}\leq\e^{2}\cdot\bigl(S_{n}(\zeta_{i})^{*}S_{n}(\zeta_{j})\bigr)_{i,j=1}^{t},
\]
we see that
\[
\norm{\pi(S)z}_{\H}^{2}\leq\e^{2}\left(\pi^{(t)}\left(\bigl(S_{n}(\zeta_{i})^{*}S_{n}(\zeta_{j})\bigr)_{i,j=1}^{t}\right)\left(y_{k}\right)_{k=1}^{t},\left(y_{\ell}\right)_{\ell=1}^{t}\right)_{\H\tensor\C^{t}}=\e^{2}\norm z_{\H}^{2}.
\]
Finally we have
\[
\norm{\pi(S)x}\le\norm{\pi(S)\left(x-z\right)}+\norm{\pi(S)z}\leq\e(\norm S+\norm z)\leq\e(\norm S+\norm x+\e),
\]
so $\pi(S)=0$ indeed.\end{proof}
\begin{defn}
\label{def:tame_subprod_sys}A subproduct system $X$ is called \emph{tame}
if $\I=\bigcap\ker\pi$, when $\pi$ ranges over all fully-coisometric
$C^{*}$-representations of $\T(X)$.\end{defn}
\begin{example}
If $X$ fulfills the requirements of \prettyref{thm:Q_n_in_T(X)}
for all $C^{*}$-representations of $\T(X)$, then every such representation
that factors through $\T(X)/\I$ is fully coisometric. Hence, by \prettyref{thm:I_subset_cap},
$X$ is tame. This class of subproduct systems is wide---see \prettyref{cor:subpr_sys_finite_dim_Hilbert_sp}
and \prettyref{prop:positive_matrix_subpr_sys}.
\end{example}

\begin{example}[cont.~of \prettyref{exa:SSP_infty}]
The subproduct system $X:=\mathrm{SSP}_{\infty}$ is tame, for if
$S\in\T(X)\backslash\I$ and $f\neq0$ is the corresponding element
of $C(B)$, there exists $z\in B'$ with $f(z)\neq0$ (as $B'$ is
dense), and $\rho_{z}$ gives rise to a fully-coisometric $C^{*}$-representation
$\pi$ of $\T(X)$ such that $\pi(S)\neq0$.\end{example}
\begin{conjecture}
\label{conj:tame_subpr_sys}All subproduct systems satisfying some
mild hypotheses are tame, at least if the adjective {}``fully-coisometric''
is replaced by {}``essential'' in \prettyref{def:tame_subprod_sys}.
\end{conjecture}
We conclude this section by giving a rough structure theory for the
representations of $\T(\mathrm{SSP}_{\infty})$.
\begin{example}[cont.~of \prettyref{exa:SSP_infty}]
For $X:=\mathrm{SSP}_{\infty}$, let $\pi$ be a $C^{*}$-representation
of $\T(X)$ on $\H$. As in \prettyref{pro:pi_I_is_pure} and the
preceding paragraph, decompose $\pi$ as $\pi_{\I}\oplus\pi_{\T(X)/\I}$,
and let $T$ be the covariant representation of $\pi_{\T(X)/\I}$.
We already know that $\pi_{\I}$ is pure (whether more could be said
is an open question). Since $\pi_{\T(X)/\I}$ factors through $\T(X)/\I\cong C(B)$,
we consider $\pi_{\T(X)/\I}$ as a $C^{*}$-representation of $C(B)$.
Write $\mathscr{K}$ for the closure of $\Img\widetilde{T}\widetilde{T}^{*}$
(equivalently, of $\Img\widetilde{T}$). Then $\mathscr{K}$ is the
closed span of the union of the images of $T(e_{n})$, $n\in\N$,
which, by virtue of normality, contains the images of $T(e_{n})^{*}$,
$n\in\N$. Thus, $\mathscr{K}$ is invariant for $\pi_{\T(X)/\I}$.
Decompose $\pi_{\T(X)/\I}$ as $\pi^{'}\oplus\pi^{''}$ with respect
to $\mathscr{K}$ and $\mathscr{K}^{\perp}$. By construction, the
$C^{*}$-representation $\pi^{'}$ is essential and $\pi^{''}$ satisfies
$\pi^{''}(S_{n}(\zeta))=0$ for all $n\in\N$ and $\zeta\in X(n)$.
\end{example}

\section{\label{sec:Morita-equivalence}Morita equivalence}

In this section we generalize ideas of \citep{Morita_equiv_tensor_alg}
to develop a notion of Morita equivalence for subproduct systems.
It is proved in Theorems \ref{thm:Morita_equiv_tensor}, \ref{thm:Morita_equiv_Toeplitz}
and \ref{thm:Morita_equiv_Cuntz_Pimsner} that if two subproduct systems
are equivalent in this sense, then so are their tensor, Toeplitz and
Cuntz-Pimsner algebras. In particular, the last theorem is proved
by showing that the Rieffel correspondence associated with the equivalence
of the Toeplitz algebras carries the ideal $\I$ of the first to that
of the second. This is yet another evidence of the naturality of the
definition of the Cuntz-Pimsner algebra for subproduct systems as
the quotient by $\I$. The results of this section should also be
compared to those of \citep{Abadie_Eilers_Exel}.

\subsection{Strong Morita equivalence of subproduct systems}

Our standard reference for Morita equivalence is \citep{Raeburn_Williams_Morita_CT}.
We assume that the reader has basic familiarity with \citep[\S 1-2]{Morita_equiv_tensor_alg},
part of which is summarized here for the sake of convenience.

Let $\Aa,\B$ be $C^{*}$-algebras, and suppose that they are Morita
equivalent via an imprimitivity bimodule $\MM$. We denote by $\widetilde{\MM}$
the opposite (dual) bimodule, and recall that the maps $m_{\Aa}:\MM\tensor_{\B}\widetilde{\MM}\to\Aa$,
$m_{\B}:\widetilde{\MM}\tensor_{\Aa}\MM\to\B$ given by $x\tensor\widetilde{y}\mapsto\leftsub{\Aa}{\langle x,y\rangle}$
and $\widetilde{x}\tensor y\mapsto\left\langle x,y\right\rangle _{\B}$,
respectively, are correspondence isomorphisms.
\begin{defn}[{\citep[Definition 2.1]{Morita_equiv_tensor_alg}}]
\label{def:correspondence_ME}Let $E,F$ be $C^{*}$-correspondences
over $\Aa,\B$, respectively. If the $C^{*}$-algebras $\Aa,\B$ are
Morita equivalent via an imprimitivity bimodule $\MM$, and if there
exists a correspondence isomorphism from $\MM\tensor_{\B}F$ onto
$E\tensor_{\Aa}\MM$, we say that $E$ and $F$ are strongly Morita
equivalent, and write $E\moritaeq_{\mathsf{M}}F$.\end{defn}
\begin{example}
If $\Aa$ and $\B$ are Morita equivalent $C^{*}$-algebras, then
they are (strongly) Morita equivalent as $C^{*}$-correspondences.
\end{example}
When the conditions of \prettyref{def:correspondence_ME} hold, the
isomorphism $W:\MM\tensor_{\B}F\to E\tensor_{\Aa}\MM$ of \prettyref{def:correspondence_ME}
induces an isomorphism $\widetilde{W}$ from $\widetilde{\MM}\tensor_{\Aa}E$
onto $F\tensor_{\B}\widetilde{\MM}$. Additionally, $E^{\tensor n}\moritaeq_{\mathsf{M}}F^{\tensor n}$
for each $n\in\N$, with correspondence isomorphisms $W_{n}:\MM\tensor_{\B}F^{\tensor n}\to E^{\tensor n}\tensor_{\Aa}\MM$
satisfying $W_{1}=W$ and
\begin{equation}
W_{n+m}=(I_{E^{\tensor n}}\tensor W_{m})(W_{n}\tensor I_{F^{\tensor m}}).\label{eq:W_n_m}
\end{equation}
Letting $W_{0}$ denote the natural isomorphism from $\MM\tensor_{\B}\B$
onto $\Aa\tensor_{\Aa}\MM$, this last equation actually holds for
all $n,m\in\Z_{+}$.

In the sequel , when $\Aa$ and $\B$ are Morita equivalent via $\MM$,
we let $L$ be the {}``linking $C^{*}$-algebra'' of $\Aa$ and
$\B$ (\citep{Brown_Green_Rieffel_1977}), namely
\[
L:=\begin{pmatrix}\B & \widetilde{\mathsf{M}}\\
\mathsf{M} & \Aa
\end{pmatrix},
\]
and for $E,F$ as above, we write $Z$ for the Hilbert $L$-module
\[
Z:=\begin{pmatrix}F & F\tensor_{\B}\widetilde{\mathsf{M}}\\
E\tensor_{\Aa}\mathsf{M} & E
\end{pmatrix}
\]
(see \citep[p.~121]{Morita_equiv_tensor_alg}).
\begin{prop}[{\citep[Proposition 2.6]{Morita_equiv_tensor_alg}}]
\label{prop:Morita_equiv_Z}If $E\moritaeq_{\mathsf{M}}F$ then there
is a left action of $L$ on $Z$, $\varphi_{Z}:L\to\L(Z)$, making
$Z$ an $L$-correspondence, satisfying $\clinspan\left(\varphi_{Z}(L)\left(\begin{smallmatrix}F & 0\\
0 & E
\end{smallmatrix}\right)\right)=Z$ and $\varphi_{Z}\left(\begin{smallmatrix}b & 0\\
0 & a
\end{smallmatrix}\right)\left(\begin{smallmatrix}\eta & 0\\
0 & \zeta
\end{smallmatrix}\right)=\left(\begin{smallmatrix}b\eta & 0\\
0 & a\zeta
\end{smallmatrix}\right)$ for $a\in\Aa$, $b\in\B$, $\zeta\in E$ and $\eta\in F$. Particularly,
$Z$ is essential.
\end{prop}
The complete definition of $\varphi_{Z}$ is given in \citep[p.~125]{Morita_equiv_tensor_alg}.

The following notion of Morita equivalence of subproduct systems is
natural in light of \prettyref{def:correspondence_ME}, as well as
\citep[Definition 5.10]{Skeide_2009_Unit_vec_Morita_equi_end}.
\begin{defn}
\label{def:subpr_sys_ME}Let $X,Y$ be subproduct systems over $\Aa,\B$
respectively, and write $E:=X(1),F:=Y(1)$. We say that $X$ and $Y$
are \emph{strongly Morita equivalent} with respect to $\mathsf{M}$
and denote $X\moritaeq_{\mathsf{M}}Y$ if $E\moritaeq_{\mathsf{M}}F$
in the sense of \prettyref{def:correspondence_ME}, with implementing
correspondence isomorphism $W:\mathsf{M}\tensor_{\B}F\to E\tensor_{\Aa}\mathsf{M}$
that satisfies
\begin{equation}
W_{n}(\mathsf{M}\tensor_{\B}Y(n))=X(n)\tensor_{\Aa}\mathsf{M},\label{eq:Morita_subproduct_sys_1}
\end{equation}
or, equivalently,
\begin{equation}
W_{n}(I_{\MM}\tensor p_{n}^{Y})=(p_{n}^{X}\tensor I_{\MM})W_{n},\label{eq:Morita_subproduct_sys_2}
\end{equation}
for all $n\in\N$. In particular, this implies that $X(n)\moritaeq_{\mathsf{M}}Y(n)$
(with $W_{n}$ implementing the equivalence). Depending upon the context,
we will regard $W_{n}$ as a mapping either from $\mathsf{M}\tensor_{\B}F^{\tensor n}$
to $E^{\tensor n}\tensor_{\Aa}\mathsf{M}$ or from $\mathsf{M}\tensor_{\B}Y(n)$
to $X(n)\tensor_{\Aa}\mathsf{M}$. The relation $\moritaeq$ is certainly
an equivalence relation.\end{defn}
\begin{rem}
\label{rem:producing_SME_subp_sys}If $X$ is a subproduct system
over $\Aa$, $F$ is an essential $C^{*}$-correspondence over $\B$
and $E:=X(1)\moritaeq_{\mathsf{M}}F$, then the implementing isomorphism
$W$ can be used to canonically induce a subproduct system $Y$ over
$\B$ with $Y(1)=F$ such that $X\moritaeq_{\mathsf{M}}Y$. Indeed,
let
\[
Y(n):=(m_{\B}\tensor I_{F^{\tensor n}})\bigl(\widetilde{\mathsf{M}}\tensor_{\Aa}W_{n}^{-1}(X(n)\tensor\MM)\bigr)
\]
for every $n\in\N$. Then $Y(n)$ is an orthogonally-complementable
sub-correspondence of $F^{\tensor n}$, $Y(n+m)\subseteq Y(n)\tensor Y(m)$
for all $n,m$, and \prettyref{eq:Morita_subproduct_sys_1} holds.
The details are left to the reader.
\end{rem}
In what follows we assume that the conditions of \prettyref{def:subpr_sys_ME}
are satisfied unless stated otherwise. For $n\in\N$, denote by $Z(n)$
the $L$-correspondence associated with the equivalence $X(n)\moritaeq_{\mathsf{M}}Y(n)$
on account of \prettyref{prop:Morita_equiv_Z},

\[
Z(n)=\begin{pmatrix}Y(n) & Y(n)\tensor_{\B}\widetilde{\mathsf{M}}\\
X(n)\tensor_{\Aa}\mathsf{M} & X(n)
\end{pmatrix},
\]
and let $Z(0):=L$ (this makes sense as $\Aa\tensor_{\Aa}\mathsf{M}\cong\MM$
and $\B\tensor_{\B}\widetilde{\mathsf{M}}\cong\widetilde{\mathsf{M}}$).
We will require the subspace $C(n):=\left(\begin{smallmatrix}Y(n) & 0\\
0 & X(n)
\end{smallmatrix}\right)$ of $Z(n)$. If $n,m\in\N$, then from \citep[Lemmas 2.7, 2.8]{Morita_equiv_tensor_alg}
we have $X(n)\tensor X(m)\moritaeq_{\mathsf{M}}Y(n)\tensor Y(m)$,
with associated $L$-correspondence
\[
Z_{n,m}:=\begin{pmatrix}Y(n)\tensor Y(m) & Y(n)\tensor Y(m)\tensor_{\B}\widetilde{\mathsf{M}}\\
X(n)\tensor X(m)\tensor_{\Aa}\mathsf{M} & X(n)\tensor X(m)
\end{pmatrix};
\]
furthermore, there is a natural $L$-correspondence isomorphism
\[
\Psi_{n,m}:Z_{n,m}\to Z(n)\tensor_{L}Z(m),
\]
which restricts to the map from $\left(\begin{smallmatrix}Y(n)\tensor Y(m) & 0\\
0 & X(n)\tensor X(m)
\end{smallmatrix}\right)$ onto $C(n)\tensor_{L}C(m)$ given by $\left(\begin{smallmatrix}\eta_{1}\tensor\eta_{2} & 0\\
0 & \zeta_{1}\tensor\zeta_{2}
\end{smallmatrix}\right)\mapsto\left(\begin{smallmatrix}\eta_{1} & 0\\
0 & \zeta_{1}
\end{smallmatrix}\right)\tensor\left(\begin{smallmatrix}\eta_{2} & 0\\
0 & \zeta_{2}
\end{smallmatrix}\right)$.
\begin{lem}
Suppose that $X,Y$ are subproduct systems over $\Aa,\B$ respectively
with $X\moritaeq_{\mathsf{M}}Y$. Using the above-mentioned notation,
the family $Z:=\left(Z(n)\right)_{n\in\Z_{+}}$ of essential $L$-correspondences
is a subproduct system: $Z(n+m)$ embeds in $Z(n)\tensor_{L}Z(m)\cong Z_{n,m}$
as an orthogonally-complementable sub-correspondence in a canonical
fashion, and the maps $(\Psi_{n,m}\tensor I_{Z(k)})\Psi_{n+m,k}$
and $(I_{Z(n)}\tensor\Psi_{m,k})\Psi_{n,m+k}$ agree on $Z(n+m+k)$,
for all $n,m,k\in\N$.\end{lem}
\begin{proof}
First, since $X,Y$ are subproduct systems, $Z(n+m)\subseteq Z_{n,m}$
as sets. We have to check that the $L$-correspondence structure of
$Z(n+m)$ (associated with $X(n+m)\moritaeq_{\mathsf{M}}Y(n+m)$)
agrees with that of $Z_{n,m}$ (associated with $X(n)\tensor X(m)\moritaeq_{\mathsf{M}}Y(n)\tensor Y(m)$
as above). To this end, we use the three formulas in the top of \citep[p.~125]{Morita_equiv_tensor_alg}.
Given $a\in\Aa$, $b\in\B$, $x,y,z,v\in\mathsf{M}$, $\zeta_{1,}\zeta_{2}\in X(n+m)$
and $\eta_{1},\eta_{2}\in Y(n+m)$, we compute:
\[
\varphi_{Z(n+m)}\begin{pmatrix}b & 0\\
0 & a
\end{pmatrix}\begin{pmatrix}\eta_{1} & \eta_{2}\tensor\widetilde{z}\\
\zeta_{1}\tensor v & \zeta_{2}
\end{pmatrix}=\begin{pmatrix}b\eta_{1} & b\eta_{2}\tensor\widetilde{z}\\
a\zeta_{1}\tensor v & a\zeta_{2}
\end{pmatrix}=\varphi_{Z_{n.m}}\begin{pmatrix}b & 0\\
0 & a
\end{pmatrix}\begin{pmatrix}\eta_{1} & \eta_{2}\tensor\widetilde{z}\\
\zeta_{1}\tensor v & \zeta_{2}
\end{pmatrix},
\]
\[
\varphi_{Z(n+m)}\begin{pmatrix}0 & 0\\
x & 0
\end{pmatrix}\begin{pmatrix}\eta_{1} & \eta_{2}\tensor\widetilde{z}\\
\zeta_{1}\tensor v & \zeta_{2}
\end{pmatrix}=\begin{pmatrix}0 & 0\\
W_{n+m}(x\tensor\eta_{1}) & (I_{X(n+m)}\tensor m_{\Aa})(W_{n+m}\tensor I_{\widetilde{\mathsf{M}}})(x\tensor\eta_{2}\tensor\widetilde{z})
\end{pmatrix},
\]
\[
\varphi_{Z(n+m)}\begin{pmatrix}0 & \widetilde{y}\\
0 & 0
\end{pmatrix}\begin{pmatrix}\eta_{1} & \eta_{2}\tensor\widetilde{z}\\
\zeta_{1}\tensor v & \zeta_{2}
\end{pmatrix}=\begin{pmatrix}(m_{\B}\tensor I_{Y(n+m)})(I_{\widetilde{\mathsf{M}}}\tensor W_{n+m}^{-1})(\widetilde{y}\tensor\zeta_{1}\tensor v) & \widetilde{W}_{n+m}(\widetilde{y}\tensor\zeta_{2})\\
0 & 0
\end{pmatrix}.
\]
The operator
\[
W_{n,m}:=(I_{X(n)}\tensor W_{m})(W_{n}\tensor I_{Y(m)}):\mathsf{M}\tensor Y(n)\tensor Y(m)\to X(n)\tensor X(m)\tensor\mathsf{M}
\]
implementing the equivalence $X(n)\tensor X(m)\moritaeq_{\mathsf{M}}Y(n)\tensor Y(m)$
extends $W_{n+m}$ by \prettyref{eq:W_n_m}, and so $\widetilde{W}_{n,m}$
extends $\widetilde{W}_{n+m}$. Thus
\[
\varphi_{Z(n+m)}\begin{pmatrix}0 & 0\\
x & 0
\end{pmatrix}\begin{pmatrix}\eta_{1} & \eta_{2}\tensor\widetilde{z}\\
\zeta_{1}\tensor v & \zeta_{2}
\end{pmatrix}=\varphi_{Z_{n,m}}\begin{pmatrix}0 & 0\\
x & 0
\end{pmatrix}\begin{pmatrix}\eta_{1} & \eta_{2}\tensor\widetilde{z}\\
\zeta_{1}\tensor v & \zeta_{2}
\end{pmatrix}
\]
and
\[
\varphi_{Z(n+m)}\begin{pmatrix}0 & \widetilde{y}\\
0 & 0
\end{pmatrix}\begin{pmatrix}\eta_{1} & \eta_{2}\tensor\widetilde{z}\\
\zeta_{1}\tensor v & \zeta_{2}
\end{pmatrix}=\varphi_{Z_{n,m}}\begin{pmatrix}0 & \widetilde{y}\\
0 & 0
\end{pmatrix}\begin{pmatrix}\eta_{1} & \eta_{2}\tensor\widetilde{z}\\
\zeta_{1}\tensor v & \zeta_{2}
\end{pmatrix}.
\]

It is easy to show that $Z(n+m)$ is orthogonally complementable in
$Z_{n,m}$: the linear mapping $p_{n,m}^{Z}:=\left(\begin{smallmatrix}p_{n+m}^{Y} & p_{n+m}^{Y}\tensor I_{\widetilde{\mathsf{M}}}\\
p_{n+m}^{X}\tensor I_{\mathsf{M}} & p_{n+m}^{X}
\end{smallmatrix}\right)$ from $Z_{n,m}$ to itself is an (orthogonal) projection in $\L(Z_{n,m})$
(for a direct calculation shows that it is a right $L$-module map),
whose range is $Z(n+m)$. In conclusion, $(\Psi_{n,m})_{|Z(n+m)}$
is an isometric, adjointable $L$-correspondence mapping from $Z(n+m)$
to $Z_{n,m}$.

For the second part of the assertion, it follows from the construction
of $\Psi$ that $(\Psi_{n,m}\tensor I_{Z(k)})\Psi_{n+m,k}$ and $(I_{Z(n)}\tensor\Psi_{m,k})\Psi_{n,m+k}$
agree on $C(n+m+k)$. Since all the maps involved are (continuous)
$L$-correspondence maps, and $\varphi_{Z(n+m+k)}(L)C(n+m+k)$ is
total in $Z(n+m+k)$ by \prettyref{prop:Morita_equiv_Z}, we infer
that the desired equality holds.\end{proof}
\begin{cor}
Under the conditions of the last lemma, $\F_{X}\moritaeq_{\mathsf{M}}\F_{Y}$
and the associated $L$-correspondence is $\F_{Z}$.
\end{cor}

\subsection{Equivalence of the operator algebras}

Let us see how the shift operator of the subproduct system $Z$, denoted
by $S^{Z}$, acts. Let $n,m\in\Z_{+}$, $\zeta_{1},\zeta_{2}\in X(n)$,
$\eta_{1},\eta_{2}\in Y(n)$, $\varrho_{1},\varrho_{2}\in X(m)$,
$\xi_{1},\xi_{2}\in Y(m)$ and $u,v,w,z\in\mathsf{M}$ be given. By
\citep[Lemma 2.9]{Morita_equiv_tensor_alg} and similar computations
that are left to the reader,
\begin{multline}
S_{n}^{Z}\begin{pmatrix}\eta_{1} & \eta_{2}\tensor\widetilde{w}\\
\zeta_{1}\tensor v & \zeta_{2}
\end{pmatrix}\begin{pmatrix}\xi_{1} & \xi_{2}\tensor\widetilde{z}\\
\varrho_{1}\tensor u & \varrho_{2}
\end{pmatrix}\\
\begin{split} & =p_{n,m}^{Z}\Psi_{n,m}^{-1}\left[\begin{pmatrix}\eta_{1} & \eta_{2}\tensor\widetilde{w}\\
\zeta_{1}\tensor v & \zeta_{2}
\end{pmatrix}\tensor_{L}\begin{pmatrix}\xi_{1} & \xi_{2}\tensor\widetilde{z}\\
\varrho_{1}\tensor u & \varrho_{2}
\end{pmatrix}\right]\\
 & =p_{n,m}^{Z}\begin{pmatrix}\eta_{1}\tensor\xi_{1}+c(\eta_{2}\tensor\widetilde{w},W_{m}^{-1}(\varrho_{1}\tensor u)) & \eta_{1}\tensor\xi_{2}\tensor\widetilde{z}+\eta_{2}\tensor\widetilde{W}_{m}(\widetilde{w}\tensor\varrho_{2})\\
\zeta_{1}\tensor W_{m}(v\tensor\xi_{1})+\zeta_{2}\tensor\varrho_{1}\tensor u & \widetilde{c}(\zeta_{1}\tensor v,\widetilde{W}_{m}^{-1}(\xi_{2}\tensor\widetilde{z}))+\zeta_{2}\tensor\varrho_{2}
\end{pmatrix}
\end{split}
\label{eq:S_Z_formula}
\end{multline}
where $c:Y(n)\tensor\widetilde{\MM}\times\MM\tensor Y(m)\to Y(n)\tensor Y(m)$
and $\widetilde{c}:X(n)\tensor\MM\times\widetilde{\MM}\tensor X(m)\to X(n)\tensor X(m)$
are given by $(\eta\tensor\widetilde{x},y\tensor\xi)\mapsto\eta\tensor\langle x,y\rangle_{\B}\xi$
and $(\zeta\tensor x,\widetilde{y}\tensor\rho)\mapsto\zeta\tensor\leftsub{\Aa}{\langle x,y\rangle}\rho$,
respectively.

The Fock space $\F_{Z}=\begin{pmatrix}\F_{Y} & \F_{Y}\tensor_{\B}\widetilde{\mathsf{M}}\\
\F_{X}\tensor_{\Aa}\mathsf{M} & \F_{X}
\end{pmatrix}$ has the following two closed linear subspaces:
\[
\F_{Z}':=\begin{pmatrix}\F_{Y} & 0\\
\F_{X}\tensor_{\Aa}\mathsf{M} & 0
\end{pmatrix},\F_{Z}'':=\begin{pmatrix}0 & \F_{Y}\tensor_{\B}\widetilde{\mathsf{M}}\\
0 & \F_{X}
\end{pmatrix}
\]
(which are left, but not right, $L$-submodules of $\F_{Z}$). From
\prettyref{eq:S_Z_formula} it is apparent that both subspaces are
invariant under the tensor algebra $\T_{+}(Z)$. As for the adjoints,
suppose that $n,m\in\N$, $z\in Z(n)$, $c\in C(n+m)$ and $l\in L$.
Approximate $c$ by a sum of the form $\sum_{i}\left(\begin{smallmatrix}\eta_{i}^{1}\tensor\eta_{i}^{2} & 0\\
0 & \zeta_{i}^{1}\tensor\zeta_{i}^{2}
\end{smallmatrix}\right)$. Then by the construction of $\Psi_{n,m}$, $S_{n}^{Z}(z)^{*}(\varphi_{Z(n+m)}(l)c)$
can be approximated by
\[
\sum_{i}\varphi_{Z(m)}\left(\left\langle z,\varphi_{Z(n)}(l)\begin{pmatrix}\eta_{i}^{1} & 0\\
0 & \zeta_{i}^{1}
\end{pmatrix}\right\rangle \right)\begin{pmatrix}\eta_{i}^{2} & 0\\
0 & \zeta_{i}^{2}
\end{pmatrix}.
\]
We therefore deduce from \prettyref{prop:Morita_equiv_Z} that $\F_{Z}'$
and $\F_{Z}''$ are also invariant under $\T_{+}(Z)^{*}$. Consequently,
they reduce the Toeplitz algebra $\T(Z)$. For convenience, we occasionally
drop the zero columns from $\F_{Z}'$ and $\F_{Z}''$.
\begin{lem}
\label{lem:T_Z_restriction_to_F_Z_tag}The restriction mappings $T\mapsto T_{|\F_{Z}'}$
and $T\mapsto T_{|\F_{Z}''}$, from $\T(Z)$ to linear operators over
$\F_{Z}'$ and $\F_{Z}''$, respectively, are injective.\end{lem}
\begin{proof}
Let $T\in\T(Z)$ be given, and suppose that $T_{|\F_{Z}'}=0$. Fix
$\zeta\in\F_{X}$, $\eta\in\F_{Y}$ and $m,z,w\in\MM$. Write $\a:=\left(\begin{smallmatrix}\eta & 0\\
0 & 0
\end{smallmatrix}\right),\be:=\left(\begin{smallmatrix}0 & 0\\
\zeta\tensor z & 0
\end{smallmatrix}\right)\in\F_{Z}'$ and $l_{1}:=\left(\begin{smallmatrix}0 & \widetilde{m}\\
0 & 0
\end{smallmatrix}\right),l_{2}:=\left(\begin{smallmatrix}0 & \widetilde{w}\\
0 & 0
\end{smallmatrix}\right)\in L$. Then
\[
T\begin{pmatrix}0 & \eta\tensor\widetilde{m}\\
0 & 0
\end{pmatrix}=T(\a\cdot l_{1})=T(\a)\cdot l_{1}=T_{|\F_{Z}'}(\a)\cdot l_{1}=0
\]
and
\[
T\begin{pmatrix}0 & 0\\
0 & \zeta\cdot\leftsub{\Aa}{\langle z,w\rangle}
\end{pmatrix}=T(\be\cdot l_{2})=T(\be)\cdot l_{2}=T_{|\F_{Z}'}(\be)\cdot l_{2}=0.
\]
Hence $T\left(\begin{smallmatrix}0 & \eta\tensor\widetilde{m}\\
0 & \zeta\langle z,w\rangle
\end{smallmatrix}\right)=0$. Since the closed span of vectors of the form $\left(\begin{smallmatrix}0 & \eta\tensor\widetilde{m}\\
0 & \zeta\langle z,w\rangle
\end{smallmatrix}\right)$ is dense in $\F_{Z}''$, we infer that $T_{|\F_{Z}''}=0$, and all
in all, $T=0$. The proof of $T_{|\F_{Z}''}=0$ $\Rightarrow$ $T=0$
is similar.
\end{proof}
Endow $\F_{Z}'$ with a right $\B$-module structure in the obvious
manner (although as a subspace of $\F_{Z}$ it is \emph{not} a right
$L$-submodule). This makes $\F_{Z}'$ a Hilbert $C^{*}$-module,
whose $\B$-valued rigging corresponds naturally to the $L$-valued
rigging of $\F_{Z}'$ as a subset of $\F_{Z}$. If $T\in\T(Z)$, it
is easy to see that $T_{|\F_{Z}'}$ is a module map, so it belongs
to $\L(\F_{Z}')$. From \prettyref{lem:T_Z_restriction_to_F_Z_tag}
it follows that the $C^{*}$-algebras homomorphism $\T(Z)\to\L(\F_{Z}')$
given by $T\mapsto T_{|\F_{Z}'}$ is injective, so that we can identify
$\T(Z)$ with its image under this map.

Denote by $\mathbf{p}$ and $\mathbf{q}$ the projections of $\F_{Z}'$
onto $\left(\begin{smallmatrix}\F_{Y}\\
0
\end{smallmatrix}\right)$ and $\left(\begin{smallmatrix}0\\
\F_{X}\tensor_{\Aa}\mathsf{M}
\end{smallmatrix}\right)$, respectively.
\begin{thm}
\label{thm:Morita_equiv_tensor}Suppose that $X,Y$ are subproduct
systems over $\Aa,\B$ respectively with $X\moritaeq_{\mathsf{M}}Y$.
Identify $\T(Z)$ with the subalgebra of $\L(\F_{Z}')$ as above.
Then:
\begin{enumerate}
\item \label{enu:Morita_equiv_tensor_1}$\mathbf{p}\T_{+}(Z)\mathbf{p}\cong\T_{+}(Y)$
and $\mathbf{q}\T_{+}(Z)\mathbf{q}\cong\T_{+}(X)$.
\item \label{enu:Morita_equiv_tensor_2}The (non-selfadjoint) operator algebras
$\T_{+}(X)$ and $\T_{+}(Y)$ are strongly Morita equivalent in the
sense of \citep{Blecher_Muhly_Paulsen__Cat_of_oper_mod}.
\end{enumerate}
\end{thm}
\begin{lem}
\label{lem:L_to_L_tensor_I_M}If $G$ is a Hilbert $C^{*}$-module
over $\Aa$ and $\MM$ is an $\Aa$-$\B$ imprimitivity bimodule,
then the map $T\mapsto T\tensor I_{\MM}$ is an isomorphism from $\L(G)$
onto $\L(G\tensor\MM)$.
\end{lem}
The proof is exactly as that of \citep[Lemma 2.12]{Morita_equiv_tensor_alg}.
The details are omitted.
\begin{proof}[Proof of \prettyref{thm:Morita_equiv_tensor}]
\ref{enu:Morita_equiv_tensor_1} Fix $n,m\in\Z_{+}$, $\zeta_{1},\zeta_{2}\in X(n)$,
$\eta_{1},\eta_{2}\in Y(n)$ and $v,w\in\mathsf{M}$. Writing
\begin{equation}
\a:=\begin{pmatrix}\eta_{1} & \eta_{2}\tensor\widetilde{w}\\
\zeta_{1}\tensor v & \zeta_{2}
\end{pmatrix}\in Z(n)\label{eq:Morita_equiv_tensor_alpha}
\end{equation}
(remember: $Z(0)=L$) we have from \prettyref{eq:S_Z_formula} that
for $\nu\in Y(m)$,
\begin{multline}
S_{n}^{Z}(\a)\begin{pmatrix}\nu\\
0
\end{pmatrix}=p_{n,m}^{Z}\begin{pmatrix}\eta_{1}\tensor\nu & 0\\
\zeta_{1}\tensor W_{m}(v\tensor\nu) & 0
\end{pmatrix}\\
=\begin{pmatrix}p_{n+m}^{Y}(\eta_{1}\tensor\nu) & 0\\
(p_{n+m}^{X}\tensor I_{\mathsf{M}})\left(\zeta_{1}\tensor W_{m}(v\tensor\nu)\right) & 0
\end{pmatrix}=\begin{pmatrix}S_{n}^{Y}(\eta_{1})\nu & 0\\
(S_{n}^{X}(\zeta_{1})\tensor I_{\mathsf{M}})W_{m}(v\tensor\nu) & 0
\end{pmatrix},\label{eq:S_Z_upper}
\end{multline}
so that
\[
\mathbf{p}S_{n}^{Z}(\a)\mathbf{p}\begin{pmatrix}\nu\\
0
\end{pmatrix}=\begin{pmatrix}S_{n}^{Y}(\eta_{1})\nu\\
0
\end{pmatrix}.
\]

Hence $\mathbf{p}\T_{+}(Z)\mathbf{p}$ is (unitarily equivalent, and
hence) completely isometrically isomorphic to $\T_{+}(Y)$. Similarly,
for $\mu\in X(m)$ and $z\in\mathsf{M}$,
\begin{equation}
\begin{split}S_{n}^{Z}(\a)\begin{pmatrix}0\\
\mu\tensor z
\end{pmatrix} & =p_{n,m}^{Z}\begin{pmatrix}c(\eta_{2}\tensor\widetilde{w},W_{m}^{-1}(\mu\tensor z)) & 0\\
\zeta_{2}\tensor\mu\tensor z & 0
\end{pmatrix}\\
 & =\begin{pmatrix}p_{n+m}^{Y}\bigl(c(\eta_{2}\tensor\widetilde{w},W_{m}^{-1}(\mu\tensor z))\bigr) & 0\\
(p_{n+m}^{X}\tensor I_{\mathsf{M}})(\zeta_{2}\tensor\mu\tensor z) & 0
\end{pmatrix}\\
 & =\begin{pmatrix}S_{n}^{Y}(\eta_{2})(m_{\B}\tensor I_{Y(m)})(\widetilde{w}\tensor W_{m}^{-1}(\mu\tensor z)) & 0\\
(S_{n}^{X}(\zeta_{2})\tensor I_{\mathsf{M}})(\mu\tensor z) & 0
\end{pmatrix},
\end{split}
\label{eq:S_Z_lower}
\end{equation}
thus
\[
\mathbf{q}S_{n}^{Z}(\a)\mathbf{q}\begin{pmatrix}0\\
\mu\tensor z
\end{pmatrix}=\begin{pmatrix}0\\
(S_{n}^{X}(\zeta_{2})\tensor I_{\mathsf{M}})(\mu\tensor z)
\end{pmatrix}.
\]
The map $T\mapsto T\tensor I_{\mathsf{M}}$ from $\L(\F_{X})$ to
$\L(\F_{X}\tensor_{\B}\mathsf{M})$ is a $C^{*}$-isomorphism by \prettyref{lem:L_to_L_tensor_I_M},
and therefore $\mathbf{q}\T_{+}(Z)\mathbf{q}$ is completely isometrically
isomorphic to $\T_{+}(X)$.

\ref{enu:Morita_equiv_tensor_2} We follow the proof of \citep[Theorem 3.2, (3)]{Morita_equiv_tensor_alg}
to show that
\[
(\mathbf{p}\T_{+}(Z)\mathbf{p},\mathbf{q}\T_{+}(Z)\mathbf{q},\mathbf{p}\T_{+}(Z)\mathbf{q},\mathbf{q}\T_{+}(Z)\mathbf{p})
\]
is a Morita context with the actions $(\mathbf{p}S_{1}\mathbf{q},\mathbf{q}S_{2}\mathbf{p}):=\mathbf{p}S_{1}\mathbf{q}S_{2}\mathbf{p}$
and $[\mathbf{q}S_{1}\mathbf{p},\mathbf{p}S_{2}\mathbf{q}]:=\mathbf{q}S_{1}\mathbf{p}S_{2}\mathbf{q}$
($S_{1},S_{2}\in\T_{+}(Z)$). The foregoing implies that $\mathbf{p},\mathbf{q}$
belong to the multiplier algebra $M(\varphi_{\infty}(L))$ and that
$\mathbf{p}\varphi_{\infty}(L)\mathbf{p}$ and $\mathbf{q}\varphi_{\infty}(L)\mathbf{q}$
are naturally isomorphic to $\varphi_{\infty}(\B)$ and $\varphi_{\infty}(\Aa)\tensor I_{\MM}$,
respectively. If $l:=\left(\begin{smallmatrix}b & \widetilde{y}\\
x & a
\end{smallmatrix}\right)\in L$, then $\mathbf{p}\varphi_{\infty}(l)\mathbf{q}\varphi_{\infty}(l)\mathbf{p}$
and $\mathbf{q}\varphi_{\infty}(l)\mathbf{p}\varphi_{\infty}(l)\mathbf{q}$
{}``equal'' $\varphi_{\infty}(\langle y,x\rangle_{\B})$ and $\varphi_{\infty}(\leftsub{\Aa}{\langle x,y\rangle})\tensor I_{\MM}$,
respectively. Consequently, the $C^{*}$-algebras $\varphi_{\infty}(\B)$
and $\varphi_{\infty}(\Aa)\tensor I_{\MM}$ are strongly Morita equivalent
through the imprimitivity bimodule $\mathbf{p}\varphi_{\infty}(L)\mathbf{q}$.
From \citep[Theorem 6.1]{Blecher_Muhly_Paulsen__Cat_of_oper_mod}
this implies that $(\varphi_{\infty}(\B),\varphi_{\infty}(\Aa)\tensor I_{\MM},\mathbf{p}\varphi_{\infty}(L)\mathbf{q},\mathbf{q}\varphi_{\infty}(L)\mathbf{p})$
is a Morita context. We omit the rest of the details, which are identical
to those of \citep{Morita_equiv_tensor_alg}.
\end{proof}
A straightforward computation using \prettyref{eq:S_Z_upper} and
\prettyref{eq:S_Z_lower} shows that $S_{n}^{Z}(\a)\mathbf{p},\mathbf{p}S_{n}^{Z}(\a)\in\T(Z)$
for all $n\in\Z_{+}$ and $\a\in Z(n)$. Hence $\mathbf{p},\mathbf{q}\in M(\T(Z))$
($\subseteq\L(\F_{Z}')$).
\begin{thm}
\label{thm:Morita_equiv_Toeplitz}Suppose that $X,Y$ are subproduct
systems over $\Aa,\B$ respectively with $X\moritaeq_{\mathsf{M}}Y$.
Identify $\T(Z)$ with the subalgebra of $\L(\F_{Z}')$ as above.
Then:
\begin{enumerate}
\item \label{enu:Morita_equiv_Toeplitz_1}$\mathbf{p}\T(Z)\mathbf{p}\cong\T(Y)$
and $\mathbf{q}\T(Z)\mathbf{q}\cong\T(X)$.
\item \label{enu:Morita_equiv_Toeplitz_2}$\T(X)\moritaeq\T(Y)$.
\end{enumerate}
\end{thm}
We shall require three technical lemmas.
\begin{lem}
\label{lem:S_and_W__X_and_Y}\mbox{}
\begin{enumerate}
\item \label{enu:S_and_W__X_and_Y_1}Let $n,k\in\N$. For all $w,w'\in\MM$,
$\eta\in Y(n)$ and $\epsilon\in Y(k)$, the operators in $\L(\F_{X}\tensor\MM)$,
defined on $X(m)\tensor\MM$, $m\in\Z_{+}$, by the formulas
\begin{equation}
W_{n+m}\left(w'\tensor\left[S_{n}^{Y}(\eta)(m_{\B}\tensor I_{Y(m)})(\widetilde{w}\tensor W_{m}^{-1}(\cdot)\right]\right)\label{eq:S_and_W__X_and_Y_1__S}
\end{equation}
and
\begin{equation}
W_{m-k}\left(w'\tensor\left[S_{k}^{Y}(\epsilon)^{*}(m_{\B}\tensor I_{Y(m)})(\widetilde{w}\tensor W_{m}^{-1}(\cdot)\right]\right)\label{eq:S_and_W__X_and_Y_1__S_star}
\end{equation}
(if $m\geq k$, otherwise $0$) can be written as $S_{n}^{X}(\zeta)\tensor I_{\MM}$
and $S_{k}^{X}(\theta)^{*}\tensor I_{\MM}$, respectively, for suitable
$\zeta\in X(n)$ and $\theta\in X(k)$.
\item \label{enu:S_and_W__X_and_Y_2}Let $n,k\in\N$. For all $v,v'\in\MM$,
$\zeta\in X(n)$ and $\xi\in X(k)$, the operators in $\L(\F_{Y})$,
defined on $Y(m)$, $m\in\Z_{+}$, by the formulas
\begin{equation}
(m_{\B}\tensor I_{Y(n+m)})\bigl(\widetilde{v'}\tensor\bigl[W_{n+m}^{-1}(S_{n}^{X}(\zeta)\tensor I_{\MM})W_{m}(v\tensor\cdot)\bigr]\bigr)\label{eq:S_and_W__X_and_Y_2__S}
\end{equation}
 and
\[
(m_{\B}\tensor I_{Y(m-k)})\bigl(\widetilde{v'}\tensor\bigl[W_{m-k}^{-1}(S_{k}^{X}(\xi)^{*}\tensor I_{\MM})W_{m}(v\tensor\cdot)\bigr]\bigr)
\]
can be written as $S_{n}^{Y}(\eta)$ and $S_{k}^{Y}(\varrho)^{*}$,
respectively, for suitable $\eta\in Y(n)$ and $\varrho\in Y(k)$.
\end{enumerate}
\end{lem}
\begin{proof}
The proofs of \ref{enu:S_and_W__X_and_Y_1} and \ref{enu:S_and_W__X_and_Y_2}
are similar, so we give details only for the former. To prove the
first part, fix $m\in\Z_{+}$, $\mu\in X(m)$ and $z\in\MM$. Approximate
$W_{m}^{-1}(\mu\tensor z)$ as the finite sum $\sum_{i}z_{i}\tensor\rho_{i}$
($z_{i}\in\MM$ and $\rho_{i}\in Y(m)$ for all $i$). Then
\[
w'\tensor\left[S_{n}^{Y}(\eta)(m_{\B}\tensor I_{Y(m)})(\widetilde{w}\tensor W_{m}^{-1}(\mu\tensor z)\right]
\]
can be approximated by
\[
(I_{\MM}\tensor p_{n+m}^{Y})\bigl(\sum_{i}w'\tensor\eta\tensor\langle w,z_{i}\rangle_{\B}\rho_{i}\bigr).
\]
Approximate $W_{n}(w'\tensor\eta)$ as the finite sum $\sum_{j}\xi_{j}\tensor x_{j}$
($\xi_{j}\in X(n)$, $x_{j}\in\MM$ for all $j$). Using \prettyref{eq:W_n_m}
and \prettyref{eq:Morita_subproduct_sys_2},
\[
\begin{split}\prettyref{eq:S_and_W__X_and_Y_1__S} & \sim\sum_{i}(p_{n+m}^{X}\tensor I_{\MM})(I_{X(n)}\tensor W_{m})\bigl(W_{n}(w'\tensor\eta)\tensor\langle w,z_{i}\rangle_{\B}\rho_{i}\bigr)\\
 & \sim\sum_{i}\sum_{j}(p_{n+m}^{X}\tensor I_{\MM})\bigl(\xi_{j}\tensor W_{m}(x_{j}\tensor\langle w,z_{i}\rangle_{\B}\rho_{i})\bigr)\\
 & =\sum_{j}(p_{n+m}^{X}\tensor I_{\MM})\bigl(\xi_{j}\tensor\leftsub{\Aa}{\langle x_{j},w\rangle}\cdot W_{m}(\sum_{i}z_{i}\tensor\rho_{i})\bigr)\\
 & \sim\sum_{j}(p_{n+m}^{X}\tensor I_{\MM})\bigl(\xi_{j}\leftsub{\Aa}{\langle x_{j},w\rangle}\tensor\mu\tensor z\bigr)\\
 & =\bigl(S_{n}^{X}(\sum_{j}\xi_{j}\leftsub{\Aa}{\langle x_{j},w\rangle})\tensor I_{\MM}\bigr)(\mu\tensor z).
\end{split}
\]
The assertion is therefore true for $\zeta:=(I_{X(n)}\tensor m_{\Aa})(W_{n}(w'\tensor\eta)\tensor\widetilde{w})$.

For the second part, fix $m\ge k$, $\mu\in X(m)$ and $z\in\MM$.
Approximate $W_{m}^{-1}(\mu\tensor z)$ as the finite sum $\sum_{i}z_{i}\tensor\rho_{i}^{(1)}\tensor\rho_{i}^{(2)}$
($z_{i}\in\MM$, $\rho_{i}^{(1)}\in Y(k)$ and $\rho_{i}^{(2)}\in Y(m-k)$
for all $i$), and $W_{k}(w\tensor\epsilon)$ as the finite sum $\sum_{j}\xi_{j}\tensor x_{j}$
($\xi_{j}\in X(k)$, $x_{j}\in\MM$ for all $j$). Then since $W_{k}$
is unitary,
\[
\begin{split}\prettyref{eq:S_and_W__X_and_Y_1__S_star} & \sim\sum_{i}W_{m-k}\left(w'\tensor\left\langle \epsilon,\langle w,z_{i}\rangle_{\B}\rho_{i}^{(1)}\right\rangle \rho_{i}^{(2)}\right)\\
 & =\sum_{i}W_{m-k}\left(w'\tensor\left\langle w\tensor\epsilon,z_{i}\tensor\rho_{i}^{(1)}\right\rangle \rho_{i}^{(2)}\right)\\
 & =\sum_{i}W_{m-k}\left(w'\tensor\left\langle W_{k}(w\tensor\epsilon),W_{k}(z_{i}\tensor\rho_{i}^{(1)})\right\rangle \rho_{i}^{(2)}\right)\\
 & \sim\sum_{i}\sum_{j}W_{m-k}\left(w'\left\langle \xi_{j}\tensor x_{j},W_{k}(z_{i}\tensor\rho_{i}^{(1)})\right\rangle \tensor\rho_{i}^{(2)}\right).
\end{split}
\]
It is easy to show that $w'\left\langle \xi\tensor x,\Theta\right\rangle =(S_{k}^{X}(\xi\cdot\leftsub{\Aa}{\langle x,w'\rangle})^{*}\tensor I_{\MM})\Theta$
for all $\xi\in X(k)$ and $\Theta\in X(k)\tensor\MM$. Thus, from
\prettyref{eq:W_n_m},
\[
\begin{split}\prettyref{eq:S_and_W__X_and_Y_1__S_star} & \sim\sum_{i}W_{m-k}\left(\left[\bigl(S_{k}^{X}(\sum_{j}\xi_{j}\cdot\leftsub{\Aa}{\langle x_{j},w'\rangle})^{*}\tensor I_{\MM}\bigr)W_{k}(z_{i}\tensor\rho_{i}^{(1)})\right]\tensor\rho_{i}^{(2)}\right)\\
 & =\left(S_{k}^{X}(\sum_{j}\xi_{j}\cdot\leftsub{\Aa}{\langle x_{j},w'\rangle})^{*}\tensor I_{\MM}\right)\sum_{i}(I_{X(k)}\tensor W_{m-k})\left(W_{k}(z_{i}\tensor\rho_{i}^{(1)})\tensor\rho_{i}^{(2)}\right)\\
 & \sim\left(S_{k}^{X}(\sum_{j}\xi_{j}\cdot\leftsub{\Aa}{\langle x_{j},w'\rangle})^{*}\tensor I_{\MM}\right)(\mu\tensor z).
\end{split}
\]
Hence $\theta:=(I_{X(k)}\tensor m_{\Aa})(W_{k}(w\tensor\epsilon)\tensor\widetilde{w'})$
fits.\end{proof}
\begin{lem}
\label{lem:S_and_W__X_and_Y_phase_2}Let $n\in\Z_{+}$. Then
\begin{enumerate}
\item for $\kappa\in Y(n)$ and $w\in\MM$, the operator in $\L(\F_{X}\tensor\MM,\F_{Y})$
defined by
\begin{equation}
X(m)\tensor\MM\ni\mu\tensor z\mapsto S_{n}^{Y}(\kappa)(m_{\B}\tensor I_{Y(m)})(\widetilde{w}\tensor W_{m}^{-1}(\mu\tensor z));\label{eq:S_and_W__X_and_Y_phase_2__1}
\end{equation}

\item and for $\varsigma\in X(n)$ and $v\in\MM$, the operator in $\L(\F_{Y},\F_{X}\tensor\MM)$
defined by
\begin{equation}
Y(m)\ni\nu\mapsto(S_{n}^{X}(\varsigma)\tensor I_{\mathsf{M}})W_{m}(v\tensor\nu);\label{eq:S_and_W__X_and_Y_phase_2__2}
\end{equation}

\end{enumerate}
belong to the closed linear span of operators of the form
\begin{equation}
X(m)\tensor\MM\ni\mu\tensor z\mapsto(m_{\B}\tensor I_{Y(n+m)})\left(\widetilde{x}\tensor W_{n+m}^{-1}\left[(S_{n}^{X}(\zeta)\tensor I_{\MM})(\mu\tensor z)\right]\right)\label{eq:S_and_W__X_and_Y_phase_2__1__clean}
\end{equation}
and
\begin{equation}
Y(m)\ni\nu\mapsto W_{n+m}(y\tensor S_{n}^{Y}(\rho)\nu),\label{eq:S_and_W__X_and_Y_phase_2__2__clean}
\end{equation}
respectively, where $\zeta\in X(n)$, $\rho\in Y(n)$ and $x,y\in\MM$.

Similar assertions are valid when $S_{n}(\cdot)$ is replaced by its
adjoint.\end{lem}
\begin{proof}
Write $\mathbf{I}$ and $\mathbf{I}\mathbf{I}$ for the operators
given by \prettyref{eq:S_and_W__X_and_Y_phase_2__1} and \prettyref{eq:S_and_W__X_and_Y_phase_2__2},
respectively. In order to make the operator $\mathbf{I}$ have the
form of \prettyref{eq:S_and_W__X_and_Y_1__S}, we ought to {}``wrap''
it with $W_{n+m}(w'\tensor\cdot)$ for some $w'\in\MM$. As $Y(n)$
is essential, $S_{n}^{Y}(\kappa)=\varphi_{\infty}(b)S_{n}^{Y}(\kappa')$
for suitable $b,\kappa'$. Since $\MM$ is full as a right $\B$-module,
$\varphi_{\infty}(b)\in\L(\F_{Y})$ belongs to the closed linear span
of operators of the form
\[
Y(p)\ni\tau\mapsto(m_{\B}\tensor I_{Y(p)})\bigl(\widetilde{x}\tensor W_{p}^{-1}(W_{p}(w'\tensor\tau))\bigr),\qquad p\in\Z_{+},
\]
where $x,w'\in\MM$. Using the first part of \prettyref{lem:S_and_W__X_and_Y},
\ref{enu:S_and_W__X_and_Y_1}, one deduces that the operator $\mathbf{I}\in\L(\F_{X}\tensor\MM,\F_{Y})$
belongs to the closed linear span of operators of the form \prettyref{eq:S_and_W__X_and_Y_phase_2__1__clean},
as stated.

To convert $\mathbf{I}\mathbf{I}$ to the form of \prettyref{eq:S_and_W__X_and_Y_2__S},
we employ the fullness of $\MM$ as a left $\Aa$-module to conclude
that for $a\in\Aa$, $\varphi_{\infty}(a)\tensor I_{\MM}\in\L(\F_{X}\tensor\MM)$
belongs to the closed linear span of operators of the form $\varphi_{\infty}(\leftsub{\Aa}{\langle y,v'\rangle})\tensor I_{\MM}$.
Notice also that
\begin{equation}
\bigl(\varphi_{\infty}(\leftsub{\Aa}{\langle y,v'\rangle})\tensor I_{\MM}\bigr)\Theta=W_{p}\left[y\tensor\left((m_{\B}\tensor I_{Y(p)})(\widetilde{v'}\tensor W_{p}^{-1}\Theta)\right)\right]\label{eq:W_and_inverse_and_m}
\end{equation}
for each $\Theta\in X(p)\tensor\MM$. As above, using the first part
of \prettyref{lem:S_and_W__X_and_Y}, \ref{enu:S_and_W__X_and_Y_2}
and that $X(n)$ is essential, the operator $\mathbf{I}\mathbf{I}\in\L(\F_{Y},\F_{X}\tensor\MM)$
belongs to the closed linear span of operators of the form \prettyref{eq:S_and_W__X_and_Y_phase_2__2__clean}.

The proof for $S_{n}(\cdot)^{*}$ goes along the lines of the preceding
one, using the second parts of \prettyref{lem:S_and_W__X_and_Y},
\ref{enu:S_and_W__X_and_Y_1} and \ref{enu:S_and_W__X_and_Y_2}.\end{proof}
\begin{lem}
\label{lem:S_Z_star}Let $k\in\N$, $\zeta_{1},\zeta_{2}\in X(k)$,
$\eta_{1},\eta_{2}\in Y(k)$ and $v,w\in\MM$. Write $\be:=\left(\begin{smallmatrix}\eta_{1} & \eta_{2}\tensor\widetilde{w}\\
\zeta_{1}\tensor v & \zeta_{2}
\end{smallmatrix}\right)\in Z(k)$. Then for $m\geq k$, $S_{k}^{Z}(\be)^{*}$ maps $\left(\begin{smallmatrix}\nu\\
\mu\tensor z
\end{smallmatrix}\right)\in Z(m)$ to
\[
\begin{pmatrix}S_{k}^{Y}(\eta_{1})^{*}\nu\\
(S_{k}^{X}(\zeta_{2})^{*}\tensor I_{\MM})(\mu\tensor z)
\end{pmatrix}+\begin{pmatrix}(m_{\B}\tensor I_{Y(m-k)})\left(\widetilde{v}\tensor W_{m-k}^{-1}\left[(S_{k}^{X}(\zeta_{1})^{*}\tensor I_{\MM})(\mu\tensor z)\right]\right)\\
W_{m-k}(w\tensor S_{k}^{Y}(\eta_{2})^{*}\nu)
\end{pmatrix}.
\]
\end{lem}
\begin{proof}
Write $S:=S_{k}^{X_{L}}(\be)$ for the shift in the full Fock space,
and remember that $S_{k}^{Z}(\be)^{*}=(S^{*})_{|\F_{Z}}$. Fix $m\geq k$.
If $\nu_{1}\in Y(k)$ and $\nu_{2}\in Y(m-k)$, then for $\nu=\nu_{1}\tensor\nu_{2}$,
\[
\begin{split}S^{*}\begin{pmatrix}\nu\\
0
\end{pmatrix} & =\varphi_{Z(m-k)}\left(\left\langle \begin{pmatrix}\eta_{1} & \eta_{2}\tensor\widetilde{w}\\
\zeta_{1}\tensor v & \zeta_{2}
\end{pmatrix},\begin{pmatrix}\nu_{1} & 0\\
0 & 0
\end{pmatrix}\right\rangle \right)\begin{pmatrix}\nu_{2} & 0\\
0 & 0
\end{pmatrix}\\
 & =\varphi_{Z(m-k)}\begin{pmatrix}\left\langle \eta_{1},\nu_{1}\right\rangle  & 0\\
w\left\langle \eta_{2},\nu_{1}\right\rangle  & 0
\end{pmatrix}\begin{pmatrix}\nu_{2} & 0\\
0 & 0
\end{pmatrix}\\
 & =\begin{pmatrix}\left\langle \eta_{1},\nu_{1}\right\rangle \nu_{2}\\
W_{m-k}(w\left\langle \eta_{2},\nu_{1}\right\rangle \tensor\nu_{2})
\end{pmatrix}=\begin{pmatrix}S_{k}^{Y}(\eta_{1})^{*}\nu\\
W_{m-k}(w\tensor S_{k}^{Y}(\eta_{2})^{*}\nu)
\end{pmatrix}.
\end{split}
\]
(see \citep[p.~125]{Morita_equiv_tensor_alg}). If $\mu_{1}\in X(k)$,
$\mu_{2}\in X(m-k)$ and $z\in\MM$, then $\Psi_{k,m-k}\left(\begin{smallmatrix}0 & 0\\
\mu_{1}\tensor\mu_{2}\tensor z & 0
\end{smallmatrix}\right)=\left(\begin{smallmatrix}0 & 0\\
0 & \mu_{1}
\end{smallmatrix}\right)\tensor_{L}\left(\begin{smallmatrix}0 & 0\\
\mu_{2}\tensor z & 0
\end{smallmatrix}\right)$ by \citep[Lemma 2.9]{Morita_equiv_tensor_alg}, and we obtain for
$\mu=\mu_{1}\tensor\mu_{2}$:
\[
\begin{split}S^{*}\begin{pmatrix}0\\
\mu\tensor z
\end{pmatrix} & =\varphi_{Z(m-k)}\left(\left\langle \begin{pmatrix}\eta_{1} & \eta_{2}\tensor\widetilde{w}\\
\zeta_{1}\tensor v & \zeta_{2}
\end{pmatrix},\begin{pmatrix}0 & 0\\
0 & \mu_{1}
\end{pmatrix}\right\rangle \right)\begin{pmatrix}0 & 0\\
\mu_{2}\tensor z & 0
\end{pmatrix}\\
 & =\varphi_{Z(m-k)}\begin{pmatrix}0 & \widetilde{v}\left\langle \zeta_{1},\mu_{1}\right\rangle \\
0 & \left\langle \zeta_{2},\mu_{1}\right\rangle
\end{pmatrix}\begin{pmatrix}0 & 0\\
\mu_{2}\tensor z & 0
\end{pmatrix}\\
 & =\begin{pmatrix}(m_{\B}\tensor I_{Y(m-k)})\left(\widetilde{v}\tensor W_{m-k}^{-1}\left[(S_{k}^{X}(\zeta_{1})^{*}\tensor I_{\MM})(\mu\tensor z)\right]\right)\\
(S_{k}^{X}(\zeta_{2})^{*}\tensor I_{\MM})(\mu\tensor z)
\end{pmatrix}.\qedhere
\end{split}
\]

\end{proof}

\begin{proof}[Proof of \prettyref{thm:Morita_equiv_Toeplitz}]
\ref{enu:Morita_equiv_Toeplitz_1} We first claim that every monomial
$S\in\T(Z)$ (say, of degree $t$) belongs to the closed linear span
of operators of the form
\begin{multline*}
Z(m)\ni\begin{pmatrix}\nu\\
\mu\tensor z
\end{pmatrix}\mapsto\\
\begin{pmatrix}T_{1}^{Y}\nu\\
(T_{1}^{X}\tensor I_{\MM})(\mu\tensor z)
\end{pmatrix}+\begin{pmatrix}(m_{\B}\tensor I_{Y(t+m)})\left(\widetilde{x'}\tensor W_{t+m}^{-1}\left[(T_{2}^{X}\tensor I_{\MM})(\mu\tensor z)\right]\right)\\
W_{t+m}(y'\tensor T_{2}^{Y}\nu)
\end{pmatrix}
\end{multline*}
for some monomials $T_{i}^{X}\in\T(X)$, $T_{i}^{Y}\in\T(Y)$ ($i=1,2$)
of degree $t$ and $x',y'\in\MM$. Indeed, suppose that $S$ is as
above. Given $n\in\Z_{+}$ and $\a\in Z(n)$, on account of \prettyref{eq:S_Z_upper},
\prettyref{eq:S_Z_lower} and \prettyref{lem:S_and_W__X_and_Y_phase_2}
we may assume that $S_{n}^{Z}(\a)$ maps $\left(\begin{smallmatrix}\nu\\
\mu\tensor z
\end{smallmatrix}\right)\in Z(m)$ to
\[
\begin{pmatrix}S_{n}^{Y}(\eta)\nu\\
(S_{n}^{X}(\xi)\tensor I_{\MM})(\mu\tensor z)
\end{pmatrix}+\begin{pmatrix}(m_{\B}\tensor I_{Y(n+m)})\left(\widetilde{x}\tensor W_{n+m}^{-1}\left[(S_{n}^{X}(\zeta)\tensor I_{\MM})(\mu\tensor z)\right]\right)\\
W_{n+m}(y\tensor S_{n}^{Y}(\rho)\nu)
\end{pmatrix}
\]
for some $\zeta,\xi\in X(n)$, $\eta,\rho\in Y(n)$ and $x,y\in\MM$.
Consequently, $S_{n}^{Z}(\a)S$ maps $\left(\begin{smallmatrix}\nu\\
\mu\tensor z
\end{smallmatrix}\right)\in Z(m)$ to
\begin{multline*}
\begin{pmatrix}S_{n}^{Y}(\eta)\left\{ T_{1}^{Y}\nu+(m_{\B}\tensor I_{Y(t+m)})\left(\widetilde{x'}\tensor W_{t+m}^{-1}\left[(T_{2}^{X}\tensor I_{\MM})(\mu\tensor z)\right]\right)\right\} \\
(S_{n}^{X}(\xi)\tensor I_{\MM})\left\{ (T_{1}^{X}\tensor I_{\MM})(\mu\tensor z)+W_{t+m}(y'\tensor T_{2}^{Y}\nu)\right\}
\end{pmatrix}\\
+\begin{pmatrix}(m_{\B}\tensor I_{Y(n+t+m)})\left(\widetilde{x}\tensor W_{n+t+m}^{-1}\left[(S_{n}^{X}(\zeta)\tensor I_{\MM})\left\{ (T_{1}^{X}\tensor I_{\MM})(\mu\tensor z)+W_{t+m}(y'\tensor T_{2}^{Y}\nu)\right\} \right]\right)\\
W_{n+t+m}(y\tensor S_{n}^{Y}(\rho)\left\{ T_{1}^{Y}\nu+(m_{\B}\tensor I_{Y(t+m)})\left(\widetilde{x'}\tensor W_{t+m}^{-1}\left[(T_{2}^{X}\tensor I_{\MM})(\mu\tensor z)\right]\right)\right\} )
\end{pmatrix}.
\end{multline*}
Utilizing \prettyref{lem:S_and_W__X_and_Y_phase_2} once again as
well as \prettyref{eq:W_and_inverse_and_m} on this last expression
yields the desired form.

We now do the same computation for the adjoints. Using \prettyref{lem:S_Z_star}
and its notation, $S_{k}^{Z}(\be)^{*}S$ maps $\left(\begin{smallmatrix}\nu\\
\mu\tensor z
\end{smallmatrix}\right)$ (when $m\ge k$) to
\begin{multline*}
\begin{pmatrix}S_{k}^{Y}(\eta_{1})^{*}\left\{ T_{1}^{Y}\nu+(m_{\B}\tensor I_{Y(t+m)})\left(\widetilde{x'}\tensor W_{t+m}^{-1}\left[(T_{2}^{X}\tensor I_{\MM})(\mu\tensor z)\right]\right)\right\} \\
(S_{k}^{X}(\zeta_{2})^{*}\tensor I_{\MM})\left\{ (T_{1}^{X}\tensor I_{\MM})(\mu\tensor z)+W_{t+m}(y'\tensor T_{2}^{Y}\nu)\right\}
\end{pmatrix}\\
+\begin{pmatrix}(m_{\B}\tensor I_{Y(t+m-k)})\left(\widetilde{v}\tensor W_{t+m-k}^{-1}\left[(S_{k}^{X}(\zeta_{1})^{*}\tensor I_{\MM})\left\{ (T_{1}^{X}\tensor I_{\MM})(\mu\tensor z)+W_{t+m}(y'\tensor T_{2}^{Y}\nu)\right\} \right]\right)\\
W_{t+m-k}(w\tensor S_{k}^{Y}(\eta_{2})^{*}\left\{ T_{1}^{Y}\nu+(m_{\B}\tensor I_{Y(t+m)})\left(\widetilde{x'}\tensor W_{t+m}^{-1}\left[(T_{2}^{X}\tensor I_{\MM})(\mu\tensor z)\right]\right)\right\} )
\end{pmatrix}.
\end{multline*}
The claim is established by appealing to \prettyref{lem:S_and_W__X_and_Y_phase_2}
and \prettyref{eq:W_and_inverse_and_m} once more.

The rest of the proof is now simple. It follows from the claim that
for every $S\in\T(Z)$ correspond $S^{X}\in\T(X)$ and $S^{Y}\in\T(Y)$
so that $\mathbf{p}S\mathbf{p}$ maps $\left(\begin{smallmatrix}\nu\\
0
\end{smallmatrix}\right)$ to $\left(\begin{smallmatrix}S^{Y}\nu\\
0
\end{smallmatrix}\right)$ and $\mathbf{q}S\mathbf{q}$ maps $\left(\begin{smallmatrix}0\\
\mu\tensor z
\end{smallmatrix}\right)$ to $\left(\begin{smallmatrix}0\\
(S^{X}\tensor I_{\MM})(\mu\tensor z)
\end{smallmatrix}\right)$. As a result, $\mathbf{p}\T(Z)\mathbf{p}$ and $\mathbf{q}\T(Z)\mathbf{q}$
are unitarily equivalent to subalgebras of $\T(Y)$ and $\T(X)\tensor I_{\MM}$,
respectively. The converse {}``inclusion'' is also true. For instance,
if $n_{1},\ldots,n_{t},m_{1},\ldots,m_{t}\in\Z_{+}$ and $\eta_{i}\in Y(n_{i})$,
$\omega_{i}\in Y(m_{i})$ for all $1\leq i\leq t$, set $\alpha_{i}:=\left(\begin{smallmatrix}\eta_{i} & 0\\
0 & 0
\end{smallmatrix}\right)$, $\be_{i}:=\left(\begin{smallmatrix}\omega_{i} & 0\\
0 & 0
\end{smallmatrix}\right)$; then $\mathbf{p}\bigl(\prod_{i=1}^{t}S_{n_{i}}^{Z}(\a_{i})^{*}S_{m_{i}}^{Z}(\be_{i})\bigr)\mathbf{p}$
{}``equals'' $\prod_{i=1}^{t}S_{n_{i}}^{Y}(\eta_{i})^{*}S_{m_{i}}^{Y}(\omega_{i})$.
This completes the proof by \prettyref{lem:L_to_L_tensor_I_M}.

\ref{enu:Morita_equiv_Toeplitz_2} We will show that $\mathbf{p}\T(Z)\mathbf{q}$
is a $\mathbf{p}\T(Z)\mathbf{p}$-$\mathbf{q}\T(Z)\mathbf{q}$ imprimitivity
bimodule, which, by the foregoing, is all we need. To this end, we
merely have to verify that $\mathbf{p}$ and $\mathbf{q}$ are full.
But we saw in the proof of \prettyref{thm:Morita_equiv_tensor}, \prettyref{enu:Morita_equiv_tensor_2},
that $\clinspan\mathbf{p}\varphi_{\infty}(L)\mathbf{q}\varphi_{\infty}(L)\mathbf{p}$
and $\clinspan\mathbf{q}\varphi_{\infty}(L)\mathbf{p}\varphi_{\infty}(L)\mathbf{q}$
{}``contain'' $\varphi_{\infty}(\B)$ and $\varphi_{\infty}(\Aa)\tensor I_{\MM}$,
respectively; and the latter sets contain approximate identities for
$\T(Y)$ and $\T(X)\tensor I_{\MM}$, respectively. Thus $\clinspan\mathbf{p}\T(Z)\mathbf{q}\T(Z)\mathbf{p}=\mathbf{p}\T(Z)\mathbf{p}$
and $\clinspan\mathbf{q}\T(Z)\mathbf{p}\T(Z)\mathbf{q}=\mathbf{q}\T(Z)\mathbf{q}$.
This completes the proof.\end{proof}
\begin{thm}
\label{thm:Morita_equiv_Cuntz_Pimsner}Suppose that $X,Y$ are subproduct
systems over $\Aa,\B$ respectively with $X\moritaeq_{\mathsf{M}}Y$.
Then $\O(X)\moritaeq\O(Y)$. More specifically, if we identify $\T(Z)$
with the subalgebra of $\L(\F_{Z}')$, $\mathbf{p}\T(Z)\mathbf{p}$
with $\T(Y)$ and $\mathbf{q}\T(Z)\mathbf{q}$ with $\T(X)$, and
treat $\mathbf{p}\T(Z)\mathbf{q}$ as a $\mathbf{p}\T(Z)\mathbf{p}$-$\mathbf{q}\T(Z)\mathbf{q}$
imprimitivity bimodule, then the image of $\I_{Y}$ under the Rieffel
correspondence (\citep[Theorem 3.22]{Raeburn_Williams_Morita_CT})
is $\I_{X}$. \end{thm}
\begin{proof}
The Morita equivalence of the Cuntz-Pimsner algebras follows from
the succeeding assertion by \citep[Proposition 3.25]{Raeburn_Williams_Morita_CT}.
Recall that $\mathbf{p}\T(Z)\mathbf{p}$ and $\mathbf{q}\T(Z)\mathbf{q}$
are \emph{naturally} unitarily equivalent to $\T(Y)$ and $\T(X)\tensor I_{\MM}$,
respectively. We have to check that
\[
\mathbf{q}\T(Z)\mathbf{p}\cdot\I_{Y}\cdot\mathbf{p}\T(Z)\mathbf{q}=\I_{X}\tensor I_{\MM}
\]
(see \citep[Proposition 3.24]{Raeburn_Williams_Morita_CT}). Since
the Rieffel correspondence is a lattice isomorphism, it is sufficient
to prove that $\mathbf{q}\T(Z)\mathbf{p}\cdot\I_{Y}\cdot\mathbf{p}\T(Z)\mathbf{q}\subseteq\I_{X}\tensor I_{\MM}$
and that $\mathbf{p}\T(Z)\mathbf{q}\cdot(\I_{X}\tensor I_{\MM})\cdot\mathbf{q}\T(Z)\mathbf{p}\subseteq\I_{Y}$.
The two inclusions are proved similarly, so we show only the first.

Let $T_{1},T_{2}\in\T(Z)$ and $S\in\I_{Y}$. Assume that $T_{2}$
is a monomial of degree $m\in\Z$. For large enough $n$, the range
of $(\mathbf{p}T_{2}\mathbf{q})(Q_{n}^{X}\tensor I_{\MM})$ is contained
in $\left(\begin{smallmatrix}Y(n+m)\\
0
\end{smallmatrix}\right)$, and so
\[
(\mathbf{q}T_{1}\mathbf{p}\cdot S\cdot\mathbf{p}T_{2}\mathbf{q})(Q_{n}^{X}\tensor I_{\MM})=(\mathbf{q}T_{1}\mathbf{p}\cdot SQ_{n+m}^{Y}\cdot\mathbf{p}T_{2}\mathbf{q})(Q_{n}^{X}\tensor I_{\MM}),
\]
and the norm of this operator is dominated by
\[
\norm{T_{1}}\norm{SQ_{n+m}^{Y}}\norm{T_{2}}\xrightarrow[n\to\infty]{}0.
\]
This proves (by \prettyref{lem:L_to_L_tensor_I_M}) that $\mathbf{q}T_{1}\mathbf{p}\cdot S\cdot\mathbf{p}T_{2}\mathbf{q}\in\I_{X}\tensor I_{\MM}$.
Since the closed span of monomials of arbitrary degree in $\T(Z)$
is $\T(Z)$, we have the desired inclusion.\end{proof}
\begin{rem}
The opposite direction, namely determining whether the Morita equivalence
of the operator algebras implies the strong Morita equivalence of
the subproduct systems, is very delicate. This is evident from the
analysis of this question in the \emph{product} system case (see \citep{Morita_equiv_tensor_alg}).
We did not attempt to tackle this problem in the present paper.
\end{rem}

\subsection{Examples}

See \citep{Morita_equiv_tensor_alg} for general examples of strong
Morita equivalence of $C^{*}$-correspondences.
\begin{example}[cf.~\citep{Muhly_Solel__Morita_Transf_of_Tensor_Alg}]
\label{exa:SME_K_C_d}Take $\Aa:=\mathbb{K}$ and $\B:=\C$, and
let $\MM$ stand for the standard $\mathbb{K}$-$\C$ imprimitivity
bimodule, namely the Hilbert space $\H:=\ell_{2}$. Fix $d\in\N$.
For a Cuntz $d$-tuple of isometries $V_{1},\ldots,V_{d}$ over $\H$
write $\a$ for the endomorphism of $\mathbb{K}$ given by $\a(T):=\sum_{i=1}^{d}V_{i}TV_{i}^{*}$.
Then $_{\a}\mathbb{K}\moritaeq_{\MM}\C^{d}$ (see \prettyref{exa:alpha_M}):
indeed, $W:{}_{\a}\mathbb{K}\tensor_{\mathbb{K}}\H\to\H\tensor_{\C}\C^{d}$
given by
\[
W(T\tensor h):=\sum_{i=1}^{d}V_{i}^{*}Th\tensor e_{i}\qquad(\A T\in\mathbb{K},h\in\H).
\]
is a correspondence isomorphism. Now $W_{n}:({}_{\a}\mathbb{K})^{\tensor n}\tensor_{\mathbb{K}}\H\to\H\tensor_{\C}(\C^{d})^{\tensor n}$
satisfies
\[
W_{n}(T_{1}\tensor\cdots\tensor T_{n}\tensor h)=\sum_{i_{1},i_{2},\ldots,i_{n}=1}^{d}V_{i_{1}}^{*}T_{1}V_{i_{2}}^{*}T_{2}\cdots V_{i_{n}}^{*}T_{n}h\tensor e_{i_{1}}\tensor e_{i_{2}}\tensor\cdots\tensor e_{i_{n}}
\]
for all $T_{1},\ldots,T_{n}\in\mathbb{K}$, $h\in\H$. Upon the identification
$({}_{\a}\mathbb{K})^{\tensor n}\cong{}_{\a^{n}}\mathbb{K}$ (which
holds since $V_{1},\ldots,V_{d}$ is {}``Cuntz''), given concretely
by
\[
T_{1}\tensor\cdots\tensor T_{n}\mapsto\sum_{i_{1},i_{2},\ldots,i_{n-1}=1}^{d}V_{i_{n-1}}\cdots V_{i_{1}}T_{1}V_{i_{1}}^{*}T_{2}V_{i_{2}}^{*}\cdots V_{i_{n-1}}^{*}T_{n},
\]
 we now get
\[
W_{n}(T\tensor h)=\sum_{i_{1},i_{2},\ldots,i_{n}=1}^{d}V_{i_{1}}^{*}\cdots V_{i_{n}}^{*}Th\tensor e_{i_{1}}\tensor e_{i_{2}}\tensor\cdots\tensor e_{i_{n}}\qquad(\A T\in\mathbb{K},h\in\H).
\]
By \prettyref{rem:producing_SME_subp_sys}, defining $Y(n)$ to be
the sub-correspondence of $_{\a^{n}}\mathbb{K}$ consisting of all
elements $T$ such that $V_{i_{1}}^{*}\cdots V_{i_{n}}^{*}T=V_{i_{\sigma(1)}}^{*}\cdots V_{i_{\sigma(n)}}^{*}T$
for every $1\leq i_{1},\ldots,i_{n}\leq d$ and $\sigma\in S_{n}$
gives a subproduct system $Y$ over $\mathbb{K}$ with $Y(1)={}_{\a}\mathbb{K}$
such that $Y\moritaeq\mathrm{SSP}_{d}$. Theorems \ref{thm:Morita_equiv_tensor},
\ref{thm:Morita_equiv_Toeplitz} and \ref{thm:Morita_equiv_Cuntz_Pimsner}
now assert that $\T_{+}(Y)$, $\T(Y)$ and $\O(Y)$ are strongly Morita
equivalent to $\T_{+}(\mathrm{SSP}_{d})$, $\T(\mathrm{SSP}_{d})$
and $\O(\mathrm{SSP}_{d})\cong C(\partial B_{d})$, respectively (each
in the appropriate sense).

More generally, all subproduct systems whose fibers are finite-dimensional
Hilbert spaces {}``come from polynomial identities'' involving finitely
many variables, and vice versa (see \citep[Proposition 7.2]{Subproduct_systems_2009}).
Hence, the construction of the last paragraph can be adapted to every
such subproduct system.
\end{example}

\begin{example}
The preceding example is valid when $B(\H)$ replaces $\mathbb{K}$.
\end{example}

\begin{example}
Take $\Aa$, $\B$, $\MM$ and $\H$ as in the last example. Fix a
Cuntz sequence $\left(V_{i}\right)_{i\in\N}$ of isometries over $\H$
with $\sum_{i=1}^{\infty}V_{i}V_{i}^{*}=I$, and let $\a$ be the
endomorphism of $B(\H)$ given by $\a(T):=\sum_{i=1}^{\infty}V_{i}TV_{i}^{*}$
(all sums are in the strong operator topology). Then $_{\a}B(\H)\moritaeq_{\MM}\H$
via $W:{}_{\a}B(\H)\tensor_{B(\H)}\H\to\H\tensor_{\C}\H$ given by
\[
W(T\tensor h):=\sum_{i=1}^{\infty}V_{i}^{*}Th\tensor e_{i}\qquad(\A T\in B(\H),h\in\H).
\]
Following the lines of \prettyref{exa:SME_K_C_d} yields a concrete
construction of a subproduct system $Y=\left(Y(n)\right)_{n\in\Z_{+}}$
over $B(\H)$ such that $Y(n)$ is a sub-correspondence of $_{\a^{n}}B(\H)$
for all $n$ and $Y\moritaeq\mathrm{SSP}_{\infty}$.
\end{example}

\section{Open questions}

In this section we state a few open questions and possible future
research directions. As usual, $X$ denotes an arbitrary (faithful)
subproduct system.
\begin{enumerate}
\item \textbf{Other characterizations of $\O(X)$.} Is there a {}``strong''
universality characterization of $\O(X)$ in the spirit of the gauge-invariant
uniqueness theorem (whether or not based on \prettyref{conj:tame_subpr_sys})?
\item \textbf{Non-faithful subproduct systems. }How should the Cuntz-Pimsner
algebra be defined for non-faithful subproduct systems? Considering
the case $X(n)=\left\{ 0\right\} $ for $n\geq n_{0}$ makes it apparent
that there is no obvious answer. Especially, it is not clear whether
adapting \prettyref{thm:I_subset_cap} to this setting is feasible.
\item \textbf{Semi-split exact sequences and $K$-theory. }Let $E$ be a
$C^{*}$-correspondence. An {}``extension of scalars'' method is
employed in \citep[\S 2]{Pimsner} to produce a $C^{*}$-correspondence
$E_{\infty}$ such that $\T(E)\hookrightarrow\T(E_{\infty})$ and
$\O(E)\cong\O(E_{\infty})$ canonically, and which admits a semi-split
exact sequence involving $\T(E_{\infty})$ and $\O(E_{\infty})$ (this
is useful for obtaining a $KK$-theoretic six-term exact sequence
for $\O(E)$). Could a similar technique be utilized for subproduct
systems? What could be said about the $K$-theory of $\O(X)$ and
$\T(X)$ in other methods (cf.~\citep[\S 8]{Katsura_2004_On_C_algebras_assoc}
and \citep{Matsumoto_1998_K_theory})?
\item \textbf{The $C^{*}$-envelope of the tensor algebra $\T_{+}(X)$.
}Denote by $C_{\mathrm{env}}^{*}(\mathcal{A})$ the $C^{*}$-envelope
of an operator algebra $\mathcal{A}$. For every $C^{*}$-correspondence
$E$ we have $C_{\mathrm{env}}^{*}(\T_{+}(E))\cong\O(E)$ by \citep[Theorem 3.7]{Katsoulis_Kribs_2006_tensor_env}.
In sharp contrast, from \citep[Theorem 8.15]{Arveson_subalgebras_3}
we obtain $C_{\mathrm{env}}^{*}(\T_{+}(\mathrm{SSP}_{d}))\cong\T(\mathrm{SSP}_{d})$.
A general statement about the relation between $C_{\mathrm{env}}^{*}(\T_{+}(X))$,
$\T(X)$ and $\O(X)$ would be very desirable.
\end{enumerate}

\section*{Acknowledgments}

The author wishes to wholeheartedly thank Baruch Solel for many stimulating
discussions, insights and suggestions. He is also very grateful to
Paul S. Muhly for his comments and ideas.

\bibliographystyle{amsplain}
\bibliography{subproduct_systems_CP}

\providecommand{\bysame}{\leavevmode\hbox to3em{\hrulefill}\thinspace}
\providecommand{\MR}{\relax\ifhmode\unskip\space\fi MR }
\providecommand{\MRhref}[2]{%
  \href{http://www.ams.org/mathscinet-getitem?mr=#1}{#2}
}
\providecommand{\href}[2]{#2}
\begin{thebibliography}{10}

\bibitem{Abadie_Eilers_Exel}
B.~Abadie, S.~Eilers, and R.~Exel, \emph{Morita equivalence for crossed
  products by {Hilbert} {$C^{*}$}-bimodules}, Trans. Amer. Math. Soc.
  \textbf{350} (1998), no.~8, 3043--3054.

\bibitem{Arveson_subalgebras_3}
W.~Arveson, \emph{Subalgebras of {$C^{*}$}-algebras {III}: {Multivariable}
  operator theory}, Acta Math. \textbf{181} (1998), no.~2, 159--228.

\bibitem{Blecher_Muhly_Paulsen__Cat_of_oper_mod}
D.~P. Blecher, P.~S. Muhly, and V.~I. Paulsen, \emph{Categories of operator
  modules ({M}orita equivalence and projective modules)}, Mem. Amer. Math. Soc.
  \textbf{143} (2000), no.~681.

\bibitem{Brown_Green_Rieffel_1977}
L.~G. Brown, P.~Green, and M.~A. Rieffel, \emph{Stable isomorphism and strong
  {M}orita equivalence of {$C\sp*$}-algebras}, Pacific J. Math. \textbf{71}
  (1977), no.~2, 349--363.

\bibitem{Davidson_C_algebras_by_example}
K.~R. Davidson, \emph{{$C^*$}-algebras by example}, Fields Institute
  Monographs, vol.~6, American Mathematical Society, Providence, RI, 1996.

\bibitem{DS2}
N.~Dunford and J.~T. Schwartz, \emph{Linear operators. {P}art {II}}, Wiley
  Classics Library, John Wiley \& Sons Inc., New York, 1988.

\bibitem{Exel_Circle_Actions}
R.~Exel, \emph{Circle actions on {$C^{*}$}-algebras, partial automorphisms, and
  a generalized {Pimsner-Voiculescu} exact sequence}, J. Funct. Anal.
  \textbf{122} (1994), no.~2, 361--401.

\bibitem{Fowler_Muhly_Raeburn_2003}
N.~J. Fowler, P.~S. Muhly, and I.~Raeburn, \emph{Representations of
  {C}untz-{P}imsner algebras}, Indiana Univ. Math. J. \textbf{52} (2003),
  no.~3, 569--605.

\bibitem{Hirshberg_2005_Ess_Rep}
I.~Hirshberg, \emph{Essential representations of {$C^*$}-correspondences},
  Internat. J. Math. \textbf{16} (2005), no.~7, 765--775.

\bibitem{Katsoulis_Kribs_2006_tensor_env}
E.~G. Katsoulis and D.~W. Kribs, \emph{Tensor algebras of
  {$C^*$}-correspondences and their {$C^*$}-envelopes}, J. Funct. Anal.
  \textbf{234} (2006), no.~1, 226--233.

\bibitem{Katsura_2003_construction}
T.~Katsura, \emph{A construction of {$C^*$}-algebras from
  {$C^*$}-correspondences}, Advances in quantum dynamics ({S}outh {H}adley,
  {MA}, 2002), Contemp. Math., vol. 335, Amer. Math. Soc., Providence, RI,
  2003, pp.~173--182.

\bibitem{Katsura_2004_On_C_algebras_assoc}
\bysame, \emph{On {$C^*$}-algebras associated with {$C^*$}-correspondences}, J.
  Funct. Anal. \textbf{217} (2004), no.~2, 366--401.

\bibitem{Katsura_2007_ideal_structure_of_C_algebras}
\bysame, \emph{Ideal structure of {$C^*$}-algebras associated with
  {$C^*$}-correspondences}, Pacific J. Math. \textbf{230} (2007), no.~1,
  107--145.

\bibitem{Lance}
E.~C. Lance, \emph{Hilbert {$C^{*}$}-modules. {A} toolkit for operator
  algebraists}, London Mathematical Society Lecture Note Series, vol. 210,
  Cambridge University Press, Cambridge, 1995.

\bibitem{Matsumoto_1997}
K.~Matsumoto, \emph{On {$C^*$}-algebras associated with subshifts}, Internat.
  J. Math. \textbf{8} (1997), no.~3, 357--374.

\bibitem{Matsumoto_1998_K_theory}
\bysame, \emph{{$K$}-theory for {$C^*$}-algebras associated with subshifts},
  Math. Scand. \textbf{82} (1998), no.~2, 237--255.

\bibitem{Tensor_algebras}
P.~S. Muhly and B.~Solel, \emph{Tensor algebras over {$C^{*}$}-correspondences:
  representations, dilations, and {$C^{*}$}-envelopes}, J. Funct. Anal.
  \textbf{158} (1998), no.~2, 389--457.

\bibitem{Morita_equiv_tensor_alg}
\bysame, \emph{On the {M}orita equivalence of tensor algebras}, Proc. Lond.
  Math. Soc. (3) \textbf{81} (2000), no.~1, 113--168.

\bibitem{Muhly_Solel__Morita_Transf_of_Tensor_Alg}
\bysame, \emph{Morita transforms of tensor algebras}, New York J. Math.
  \textbf{17a} (2011), 87--100.

\bibitem{Paschke_1973}
W.~L. Paschke, \emph{Inner product modules over {$B^{*}$}-algebras}, Trans.
  Amer. Math. Soc. \textbf{182} (1973), 443--468.

\bibitem{Pimsner}
M.~V. Pimsner, \emph{A class of {$C^*$}-algebras generalizing both
  {C}untz-{K}rieger algebras and crossed products by {${\mathbb Z}$}}, Free
  probability theory ({W}aterloo, {ON}, 1995) (D.~Voiculescu, ed.), Fields
  Inst. Commun., vol.~12, Amer. Math. Soc., Providence, RI, 1997, pp.~189--212.

\bibitem{Raeburn_Williams_Morita_CT}
I.~Raeburn and D.~P. Williams, \emph{Morita equivalence and continuous-trace
  {$C^*$}-algebras}, Mathematical Surveys and Monographs, vol.~60, American
  Mathematical Society, Providence, RI, 1998.

\bibitem{Subproduct_systems_2009}
O.~M. Shalit and B.~Solel, \emph{Subproduct systems}, Doc. Math. \textbf{14}
  (2009), 801--868.

\bibitem{Sims_Yeend_2010}
A.~Sims and T.~Yeend, \emph{{$C^*$}-algebras associated to product systems of
  {H}ilbert bimodules}, J. Operator Theory \textbf{64} (2010), no.~2, 349--376.

\bibitem{Skeide_2009_Unit_vec_Morita_equi_end}
M.~Skeide, \emph{Unit vectors, {M}orita equivalence and endomorphisms}, Publ.
  Res. Inst. Math. Sci. \textbf{45} (2009), no.~2, 475--518.

\bibitem{Viselter_cov_rep_subproduct_systems}
A.~Viselter, \emph{Covariant representations of subproduct systems}, Proc.
  Lond. Math. Soc. (3) \textbf{102} (2011), no.~4, 767--800.

\end{thebibliography}

\end{document}